\documentclass[12pt,leqno]{amsart}
\usepackage{amssymb}
\usepackage{amsmath,amssymb,color}
\numberwithin{equation}{section}

\newcommand{\ep}{\varepsilon}

\newcommand{\ppp}{\partial}

\newcommand{\pppa}{\partial_t^{\alpha}}

\newcommand{\hhalf}{\frac{1}{2}}

\newcommand{\DDD}{\mathcal{D}}

\newcommand{\ddda}{d_t^{\alpha}}
\newcommand{\DDDa}{D_t^{\alpha}}

\newcommand{\R}{\mathbb{R}}

\newcommand{\N}{\mathbb{N}}

\newcommand{\www}{\widetilde}

\newcommand{\CC}{{_{0}{C^1[0,T]}}}

\newcommand{\ooo}{\overline}

\newcommand{\BLS}{\blacksquare}

\newcommand{\RRRR}{\longrightarrow}

\newcommand{\WWW}{W_{\alpha,p}(0,T)}
\newcommand{\WWWEP}{W^{\alpha-\varepsilon,p}(0,T)}
\newcommand{\LLLP}{L^p(0,T)}
\newcommand{\LLLQ}{L^q(0,T)}
\newcommand{\CCCM}{\,_{0}C^m[0,T]}
\newcommand{\WWWR}{\,_{\alpha,p}W(0,T)}
\newcommand{\WWWRO}{\,^{0}W^{\alpha,p}(0,T)}

\newcommand{\JJJ}{J^{\alpha}}
\allowdisplaybreaks

\title
[]
{Fractional derivatives and time-fractional ordinary differential equations
in $L^p$-space}

\author{
$^{1,2,3}$ M.~Yamamoto
}

\thanks{
\\
$^1$ Graduate School of Mathematical Sciences, The University
of Tokyo, Komaba, Meguro, Tokyo 153-8914, Japan \\
$^2$ Honorary Member of Academy of Romanian Scientists, 
Ilfov, nr. 3, Bucuresti, Romania \\
$^3$ Correspondence Member of Accademia Peloritana dei Pericolanti,
Palazzo Universit\`a, Piazza S. Pugliatti 1 98122 Messina, Italy\\
$^4$ Peoples' Friendship University of Russia 
(RUDN University) 6 Miklukho-Maklaya St, Moscow, 117198, Russian Federation
\\
e-mail: {\tt myama@ms.u-tokyo.ac.jp}
}

\date{}
\begin{document}
\maketitle

\baselineskip 18pt

\begin{abstract}
We define fractional derivatives $\pppa$ in Sobolev spaces based on 
$L^p(0,T)$ by an operator theory, and characterize the domain of
$\pppa$ in subspaces of the Sobolev-Slobodecki spaces $W^{\alpha,p}(0,T)$.
Moreover we define $\pppa u$ for $u\in L^p(0,T)$ in a sense of 
distribution.  Then we discuss initial value problems for 
linear fractional ordinary differential equations by means of 
such $\pppa$ and establish several results on the unique existence of 
solutions within specified classes according to the regularity of 
the coefficients and the non-homogeneous terms in the equations.
\end{abstract}

\section{Introduction}

Recently fractional differential equations have attracted great attention, 
not only by theoretical interests but also for by the necessity for 
various applications.  For example, time-fractional differential equations
are understood as reasonable model equations in various anomalous diffusion
phenomena which are frequently observed in diffusion in heterogeneous media 
around us (e.g., Metzler and Klafter \cite{MK}, Chapter 10 in 
Podlubny \cite{Po}).

For $0 < \alpha < 1$, we can formally define 
the pointwise Caputo derivative by 
$$
\ddda v(t) =
\frac{1}{\Gamma(1-\alpha)}\int^t_0 (t-s)^{-\alpha}
\frac{dv}{ds}(s) ds,                       \eqno{(1.1)}
$$
as long as the right-hand side exists.
Here and henceforth, for $\beta > 0$ by $\Gamma(\beta)$ we denote 
the gamma function: $\Gamma(\beta):= \int^{\infty}_0 e^{-t}t^{\beta-1}dt$.
As for classical treatments on fractional derivatives, 
we can refer to monographs Gorenflo, Kilbas, Mainardi and Rogosin 
\cite{GKMR}, Kilbas, Srivastava and Trujillo \cite{KST},
Podlubny \cite{Po} for example.
  
We can interpret $\ddda v$ as the $\alpha$-th derivaive of $v$, but
in (1.1), the first derivative $\frac{dv}{ds}$ appears and so
intuitively $\ddda v$ requires the existence of $\frac{dv}{dt}$ 
even for defining the lower-order derivative of order $\alpha<1$.
Moreover, a corresponding initial value problem for a time-fractional 
ordinary differential equation can be formulated as 
$$
\ddda u(t) = b(t)u(t) + f(t), \quad 0<t<T, 
                                       \eqno{(1.2)}
$$
and
$$
u(0) =a.                               \eqno{(1.3)}
$$
Since only $\alpha$-th derivative appears with $0<\alpha<1$ in (1.2), 
the initial condition (1.3) is not trivial.  This can be understood
in a simple case
where $b\equiv 0$ and $f(t) = t^{-\gamma}$ with $0<\gamma<1$.  Then
$$
u(t) = \frac{\Gamma(1-\gamma)}{\Gamma(\alpha-\gamma+1)}t^{\alpha-\gamma}
+ a
$$
satisfies (1.2) - (1.3) for $0<\gamma<\alpha<1$, but cannot satisfy 
(1.2) - (1.3) if $0<\alpha<\gamma<1$.

This observation suggests us to make adequate formulations for
both fractional derivatives and initial value problems, when we consider
not smooth functions.
There are several studies and here we refer to 
In Gorenflo, Luchko and Yamamoto \cite{GLY} and Kubica, Ryszewska and 
Yamamoto \cite{KRY}, Yamamoto \cite{Ya21}, where one constructs the
fractional derivative operator in Sobolev spaces based on 
$L^2(0,T)$ and apply it to fractional differential equations 
including partial differential equations in spatial variables.
Also see Kian \cite{Ki}, Kian and Yamamoto \cite{KiYa},
Luchko and Yamamoto \cite{LY}, Sakamoto and Yamamoto \cite{SY},
Zacher \cite{Za}.
The method in those works relies on the structure of Hilbert spaces,
that is, $p=2$.

The purpose of this article is to establish the corresponding theory based on 
$L^p(0,T)$ with $1\le p < \infty$.
\\

Throughout this article, we assume
$$
1 \le p < \infty.
$$
Henceforth let 
$L^p(0,T):= \{ v;\, \int_0^T \vert v(t)\vert^p dt < \infty\}$ with 
$1\le p < \infty$
and $W^{1,p}(0,T):= \{ v\in L^p(0,T);\, \frac{dv}{dt} \in L^p(0,T)\}$.
We define the norm by 
$$
\Vert v\Vert_{L^p(0,T)}:= \left(\int_0^T \vert v(t)\vert^p dt\right)
^{\frac{1}{p}}, \quad
\Vert v\Vert_{W^{1,p}(0,T)}:= \Vert v\Vert_{L^p(0,T)}
+ \left\Vert \frac{dv}{dt}\right\Vert_{L^p(0,T)}.
$$

For $\alpha>0$, we set
$$
J^{\alpha}v(t) = \frac{1}{\Gamma(\alpha)}\int^t_0 (t-s)^{\alpha-1}
v(s) ds, \quad 0<t<T, \, v\in L^1(0,T),      \eqno{(1.4)}
$$
which is called the Riemann-Liouville fractional integral operator, and 
can be interpreted as an $\alpha$-times integral of $v$.
By the Young inequality on the convolution (e.g., Lemma A.1 in 
\cite{KRY}), we see that 
$$
J^{\alpha}L^q(0,T) \subset L^q(0,T) \quad \mbox{with $q\ge 1$}.
                                       \eqno{(1.5)}
$$

Our definition of the fractional derivative relies on 
the inverse to $J^{\alpha}$ in $L^p(0,T)$.
Henceforth by $\pppa$ we denote such a fractional derivative in order to 
distinguish from $\ddda$.  The first step is to first clarify the 
space $J^{\alpha}L^p(0,T)$ and introduce a suitable norm.
After establishing $\pppa$, we formulate an initial value problem
and discuss the unique existence of solutions with suitable 
regularity.

In order to accomplish the task, we organize the article as follows:
\begin{itemize}
\item
\S2. Fractional derivative $\pppa$ in $L^p$-Sobolev-Slobodecki spaces 
of positive orders
\item
\S3. Fractional derivative $\pppa$ in $L^p(0,T)$
\item
\S4. Initial value problems for time-fractional ordinary differential
equations
\item
\S5. Concluding remarks
\item
\S6. Appendix: Proof of Lemma 2.4.
\end{itemize}

\section{Fractional derivative $\pppa$ in $L^p$-Sobolev-Slobodecki spaces 
of positive orders: characterization of the domain $\DDD(\pppa)$ and
norm estimates}

Henceforth by $\DDD(K)$, we denote the domain of an operator $K$ under 
consideration, and $C>0$ denotes generic constants depending on 
$\alpha, p$ but not on choices of functions $v$.

We recall that $J^{\alpha}$ is defined by (1.4) and we set 
$\DDD(J^{\alpha}) = L^p(0,T)$.
Then by (1.5) we see that 
$$
\Vert J^{\alpha}v\Vert_{L^p(0,T)} \le C\Vert v\Vert_{L^p(0,T)}
\quad \mbox{for all $v\in L^p(0,T)$}.                          \eqno{(2.1)}
$$

First we prove:
\\
{\bf Lemma 2.1.}
\\
{\it  Let $\alpha, \beta > 0$ and $\alpha, \beta \not\in \N$.
\\
(i) $\JJJ J^{\beta}v = J^{\alpha+\beta}v$ for each $v\in L^1(0,T)$.
\\
(ii) $\JJJ: \LLLP \, \RRRR\, \LLLP$ is injective for $1\le p < \infty$..
}
\\
{\bf Proof of Lemma 2.1.}
\\
(i) Exchanging the order of the integral, we have 
\begin{align*}
& J^{\alpha}(J^{\beta}v)(t)
= \frac{1}{\Gamma(\alpha)}\int^t_0 (t-s)^{\alpha-1}
J^{\beta}v(s) ds\\
=& \frac{1}{\Gamma(\alpha)}\int^t_0 (t-s)^{\alpha-1}
\frac{1}{\Gamma(\beta)}\left( \int^s_0 (s-\xi)^{\beta-1}v(\xi) d\xi
\right) ds\\
=& \frac{1}{\Gamma(\alpha)\Gamma(\beta)}\int^t_0 
\left( \int^t_{\xi} (t-s)^{\alpha-1}(s-\xi)^{\beta-1}ds \right)
v(\xi) d\xi\\
= &\frac{1}{\Gamma(\alpha)\Gamma(\beta)}\int^t_0 
\frac{\Gamma(\alpha)\Gamma(\beta)}{\Gamma(\alpha+\beta)}
(t-\xi)^{\alpha+\beta-1}v(\xi) d\xi
= (J^{\alpha+\beta}v)(t), \quad 0<t<T,
\end{align*}
which completes the proof of (i).  $\blacksquare$
\\
(ii) We assume that $v\in L^1(0,T)$ satisfies 
$$
J^{\alpha}v(t) = \frac{1}{\Gamma(\alpha)}
\int^t_0 (t-s)^{\alpha-1} v(s) ds = 0, \quad 0<t<T.
$$
We set $\alpha = (m-1) + \sigma$ with $m\in \N$ and 
$0<\sigma<1$.  We operate $J^{1-\sigma}$ to obtain
$$
\int^t_0 (t-s)^{-\sigma}\left(\int^s_0 (s-\xi)^{\alpha-1} v(\xi) d\xi
\right) ds = 0, \quad 0<t<T.
$$
Again exchanging the order of the integral and using $\alpha-\sigma
= m-1$, we obtain
\begin{align*}
& 0 = \int^t_0 \left( \int^t_{\xi} (t-s)^{-\sigma}(s-\xi)^{\alpha-1} ds
\right) v(\xi) d\xi\\
=& \frac{\Gamma(1-\sigma)\Gamma(\alpha)}{\Gamma(1-\sigma+\alpha)}
\int^t_0 (t-\xi)^{m-1} v(\xi) d\xi, \quad 0<t<T.
\end{align*}
Therefore, the $m$-times differentiation yields $v(t) = 0$ for 
$0<t<T$.  Thus the proof of Lemma 2.1 is complete.
$\blacksquare$
\\

The proof of (ii) is suggested by the classical formula 
(e.g., Gorenflo and Vessella \cite{GV}) 
$$
D_t^{\alpha}J^{\alpha}v = v, \quad v\in L^1(0,T),  
$$
where for $\alpha = (m-1) + \sigma$ with 
$m\in \N$ and $0<\sigma<1$, 
we define the Riemann-Liouville fractional derivative
$$
\DDDa v(t) = \frac{1}{\Gamma(1-\sigma)} \frac{d^m}{dt^m}
\int^t_0 (t-s)^{-\sigma} v(s) ds, \quad v\in C^{\infty}[0,T].
                                                         \eqno{(2.2)}
$$
\\

By Lemma 2.1 (ii), there exists an inverse operator $(J^{\alpha})^{-1}$
in $\LLLP$ to $J^{\alpha}$ algebraically.

We define  
$$
\WWW:= \JJJ\LLLP \quad \mbox{as a set}       \eqno{(2.3)}
$$
with the norm
$$
\Vert u\Vert_{\WWW}:= \Vert (J^{\alpha})^{-1}u\Vert_{L^p(0,T)}.
                                              \eqno{(2.4)}
$$
We note that
$$
\Vert J^{\alpha}w\Vert_{\WWW} = \Vert w\Vert_{L^p(0,T)}.
$$
Therefore, 
\\
{\bf Lemma 2.2.}
\\
{\it
Let $\alpha > 0$.  If 
$w_n \RRRR w$ in $\LLLP$, then $J^{\alpha}w_n \RRRR J^{\alpha}w$ in 
$\WWW$.
}

Moreover, we prove
\\
{\bf Lemma 2.3.}
\\
{\it
The space $\WWW$ is complete with the norm $\Vert \cdot\Vert_{\WWW}$, that is,
$\WWW$ is a Banach space.
}
\\
{\bf Proof.}
\\
Let $u_n \in \WWW$ and $\lim_{m,n\to\infty} \Vert u_n-u_m\Vert_{\WWW} = 0$.
Then for each $n\in \N$, we can find $w_n \in \LLLP$ such that 
$u_n = \JJJ w_n$.  Hence, 
$$
\Vert u_m-u_n\Vert_{\WWW} = \Vert \JJJ(w_m-w_n)\Vert_{\WWW}
= \Vert w_m-w_n\Vert_{\LLLP}\, \RRRR 0
$$
as $m,n \to \infty$.  By the completeness of $\LLLP$, there exists 
$w_0 \in \LLLP$ such that $\lim_{n\to\infty} \Vert w_n-w_0\Vert_{\LLLP}
= 0$.  Set $u_0:= \JJJ w_0$.
Lemma 2.2 yields $\JJJ w_n \RRRR \JJJ w_0$ in $\WWW$, that is,
$u_n \RRRR \JJJ w_0$ in $\WWW$, where $\JJJ w_0 \in \JJJ\LLLP
= \WWW$.  Thus the proof of Lemma 2.3 is complete.
$\blacksquare$
\\
\vspace{0.1cm}

We define the fractional derivative.
\\
{\bf Defintion 2.1.}
\\
{\it 
Let $\alpha > 0$ and $\alpha \not\in \N$.  Then
$$
\pppa:= (J^{\alpha})^{-1} \quad \DDD(\pppa) = \WWW.
$$
}

At least $\pppa$ is well-defined, but is not useful if we do not specify
the domain $\DDD(\pppa)$.   
In a special case of $p=2$, by the structure of the Hilbert space 
$L^2(0,T)$, we can make the complete characterization in terms of
$W^{\alpha,2}(0,T)$ (e.g., \cite{GLY}, \cite{KRY}, \cite{Ya21}).
However for $p\ne 2$, to the author's best knowledge, such characterization
is not known.

Now we will make characterization of $\WWW$ and for it we introduce 
function spaces.
For $m\in \N$, we set $\frac{d^0v}{dt^0}(t) = v(t)$ and 
$$
_{0}C^m[0,T]:= 
\left\{ v\in C^m[0,T];\, \frac{d^kv}{dt^k}(0) = 0\quad 
\mbox{for $k=0,1,..., m-1$}\right\}.
$$
For $\alpha = m-1 + \sigma$ with $m\in \N$ and $0<\sigma<1$ and
$1 \le p < \infty$, we define the Sobolev-Slobodecki space 
$W^{\alpha,p}(0,T)$ with the norm $\Vert \cdot\Vert_{W^{\alpha,p}(0,T)}$ 
by 
$$
W^{\alpha,p}(0,T) := \{ u\in L^p(0,T);\, 
\Vert u\Vert_{W^{\alpha,p}(0,T)} < \infty\}
$$
where 
$$
\Vert u\Vert_{W^{\alpha,p}(0,T)}:=
\left( \Vert u\Vert^p_{W^{m-1,p}(0,T)}
+ \int^T_0\int^T_0 \frac{\left\vert \frac{d^{m-1}u}{dt^{m-1}}(t)
- \frac{d^{m-1}u}{dt^{m-1}}(s)\right\vert^p}{\vert t-s\vert^{1+\sigma p}}
dtds \right)^{\frac{1}{p}}                
$$
(e.g., Adams \cite{Ad}, Grisvard \cite{Gri}).  
By $\ooo{X}^{Y}$ we denote the closure of $X\subset Y$ in a normed space 
$Y$.
Moreover
$$
_{0}W^{\alpha,p}(0,T):= \, \ooo{_{0}C^m[0,T]}^{W^{\alpha,p}(0,T)}.
$$
We set $W^{0,p}(0,T) = \, _{0}W^{\alpha,p}(0,T) = L^p(0,T)$ for uniform 
notations.
\\

We state the main result in this section.
\\
{\bf Theorem 2.1.}
\\
{\it 
Let $\alpha = m-1+\sigma$ with $m\in \N$ and $0<\sigma<1$, and
$\ep > 0$ be arbitrarily given such that $0<\ep<\alpha$.
Then
\\
(i) There exists a constant $C=C(\alpha,\ep) > 0$ such that 
$$
\Vert J^{\alpha}v\Vert_{W^{\alpha-\ep,p}(0,T)} \le C\Vert v\Vert_{L^p(0,T)}
                                                      \eqno{(2.5)}
$$
for all $v\in \LLLP$ and 
$$
\WWW \subset \, _{0}\WWWEP.
$$
\\
(ii)
We have
$$
_{0}W^{\alpha+\ep,p}(0,T) \subset \WWW 
$$
and
$$
J^{\alpha}\left( \frac{d^m}{dt^m}J^{m-\alpha}u\right)(t) = u(t) \quad 
\mbox{for each $u \in \, _{0}W^{\alpha+\ep,p}(0,T)$}.
                                                      \eqno{(2.6)}
$$
\\
(iii) The embedding $W_{\alpha+\ep,p}(0,T) \RRRR \WWW$ is compact.
}
\\

In terms of the Riemann-Liouville fractional derivative (2.2), we can 
rewrite (2.6) as
$$
J^{\alpha}D_t^{\alpha}u = u \quad
\mbox{for each $u \in \, _{0}W^{\alpha+\ep,p}(0,T)$}.
$$

We note that for $\alpha = m-1$ with $m\in \N$, we can easily prove
$$
J^{\alpha}L^p(0,T) 
= \left\{ v\in W^{m,p}(0,T);\, \frac{d^kv}{dt^k}(0) = 0\quad 
\mbox{for $k=0,1,..., m-1$}\right\}.
$$
In the case of $p=2$ and $0<\alpha<1$, it is proved in 
\cite{KRY} that $W_{\alpha,2}(0,T)
= \,\ooo{_{0}C^1[0,T]}^{W^{\alpha,2}(0,T)}$, which means that we can 
take $\ep=0$ in (i) and (ii).
However, in the case $p \ne 2$, we do not know the equivalent representation of
$\WWW$, but (i) and (ii) imply 
$$
_{0}W^{\alpha+\ep,p}(0,T) \subset \WWW \subset \, _{0}\WWWEP
$$
with arbitrarily small $\ep \in (0,\alpha)$, which means that 
$\WWW$ is interpolated between $_{0}W^{\alpha+\ep,p}(0,T)$ and
$_{0}W^{\alpha-\ep,p}(0,T)$ with arbitrarly small gap $\ep>0$.
\\
\vspace{0.1cm}
\\
{\bf Corollary 2.1.}
\\
{\it
Let $\alpha>0$ and $\beta \ge 0$.  Then
\\
(i) $J^{\alpha}:\, W_{\beta,p}(0,T) \RRRR W_{\beta,p}(0,T)$ is 
a compact operator.
\\
(ii) $J^{\alpha}W_{\beta,p}(0,T) = W_{\alpha+\beta,p}(0,T)$ and
$\Vert J^{\alpha}w\Vert_{W_{\alpha+\beta,p}(0,T)}
= \Vert w\Vert_{W_{\beta,p}(0,T)}$ for $w\in W_{\beta,p}(0,T)$.
}
\\

We rewrite Theorem 2.1 in terms of $\pppa = (J^{\alpha})^{-1}$.
\\
{\bf Theorem 2.2.}
\\
{\it
Let $p \ge 1$, $\alpha = m-1 + \sigma$ with $m\in \N$ and $0<\sigma<1$,
and $0<\ep<\alpha$.
\\
(i) $\, _{0}W^{\alpha+\ep,p}(0,T) \subset \DDD(\pppa) \, \subset \,
_{0}W^{\alpha-\ep,p}(0,T)$ topologically and algebraically.
\\
In particular, there exists a constant $C=C(\alpha,\ep) > 0$ such that 
$$
C^{-1}\Vert u\Vert_{W^{\alpha-\ep,p}(0,T)} 
\le \Vert \pppa u\Vert_{\LLLP}
\le C\Vert \pppa u\Vert_{W^{\alpha+\ep,p}(0,T)}
$$
for all $u \in \, _{0}W^{\alpha+\ep,p}(0,T)$.
\\
(ii) $\pppa W_{\alpha+\beta,p}(0,T) = W_{\beta,p}(0,T)$ for 
$\alpha > 0$ and $\beta \ge 0$.
}

The part (i) of Theorem 2.1 improves Proposition 6 in 
Carbotti and Comi \cite{CC}, and see Li and Liu \cite{LL1}, \cite{LL2}
as for results on fractional derivative in Sobolev spaces on $L^p(0,T)$.

The rest of this section is devoted to the proof of Theorem 2.1,
but the arguments in the proof are not used later.
The proof of Theorem 2.1 (i) in the case of $p=1$ and 
$0<\alpha<1$ can be found in the proof of Theorem 4.2.2
(pp.73-74) in Gorenflo and Vessella \cite{GV}, and we adjust 
their arguments in \cite{GV} to the general case $p>1$.
Moreover in K\"onig \cite{Koe}, related results are discussed in an 
operator theoretical setting, while our proof is lengthy 
but based on direct estimation.
\\
\vspace{0.1cm}
\\
{\bf Proof of Theorem 2.1 (i).}
\\
{\bf First Step: $0<\alpha<1$}.
\\
First we show that for $\delta \in \left( 0, \hhalf\right)$, 
there exists a constant $C_{\delta}>0$ such that 
$$
\vert 1-\eta^{\alpha-1}\vert \le (1-\alpha)(1-\eta)
+ C_{\delta}(1-\eta)^2 \quad \mbox{if $1-\delta \le
\eta \le 1$}.                                       \eqno{(2.7)}
$$
\\
{\bf Verification of (2.7).}
\\
Applying the Taylor theorem to $g(\eta):= \eta^{\alpha-1}$,
we obtain
$$
g(\eta) = g(1) + g'(1)(\eta-1) + \frac{g''(\zeta)}{2!}(\eta-1)^2
\quad \mbox{for $1-\delta \le \eta \le 1$},
$$
where $\zeta \in (1-\delta,1)$ is some number. Therefore
$\vert \eta^{\alpha-1} - 1 - (\alpha-1)(\eta-1)\vert 
\le C_{\delta}(\eta-1)^2$, which implies (2.7).
$\BLS$
\\ 

We set $\beta := \alpha-\ep$ with $0<\ep<\alpha$.
We recall that 
$$
\Vert v\Vert^p_{W^{\beta,p}(0,T)}
= \Vert v\Vert^p_{\LLLP}
+ \int^T_0\int^T_0 \frac{\vert v(t) - v(s)\vert^p}
{\vert t-s\vert^{1+\beta p}} dtds.
$$
We can see that it suffices to estimate for $0<s<t$.
Indeed, dividing $[0,T]^2
= \{ (t,s)\in [0,T]^2;\, 0\le s\le t \le T\}
\cup \{ (t,s)\in [0,T]^2;\, 0\le t\le s \le T\}$, we 
have 
$$
\int^T_0\int^T_0 \frac{\vert v(t) - v(s)\vert^p}
{\vert t-s\vert^{1+\beta p}} dtds
= 2\int^T_0 \left( \int^t_0 \frac{\vert v(t) - v(s)\vert^p}
{\vert t-s\vert^{1+\beta p}} ds \right) dt.
$$
We have 
\begin{align*}
& \vert \Gamma(\alpha)(J^{\alpha}v(t) - J^{\alpha}v(s))\vert
= \left\vert \int^t_0 (t-\xi)^{\alpha-1} v(\xi)d\xi
- \int^s_0 (s-\xi)^{\alpha-1} v(\xi) d\xi\right\vert\\
= & \left\vert \int^t_s (t-\xi)^{\alpha-1} v(\xi)d\xi
+ \int^s_0 ((t-\xi)^{\alpha-1} - (s-\xi)^{\alpha-1}) v(\xi) d\xi
\right\vert\\
\le & \int^t_s (t-\xi)^{\alpha-1} \vert v(\xi)\vert d\xi
+ \int^s_0 \vert (t-\xi)^{\alpha-1} - (s-\xi)^{\alpha-1}\vert 
\vert v(\xi)\vert  d\xi,
\end{align*}
and so
\begin{align*}
& \vert J^{\alpha}v(t) - J^{\alpha}v(s)\vert^p\\
\le & C\left( \int^t_s (t-\xi)^{\alpha-1} \vert v(\xi)\vert d\xi
\right)^p
+ C\left( \int^s_0 \vert (t-\xi)^{\alpha-1} - (s-\xi)^{\alpha-1}\vert 
\vert v(\xi)\vert  d\xi\right)^p.
\end{align*}
Hence
$$
\int^T_0\int^t_0
\frac{\vert J^{\alpha}v(t) - J^{\alpha}v(s)\vert^p}
{\vert t-s\vert^{1+\beta p}} dsdt
\le C\int^T_0 dt \int^t_0 ds (t-s)^{-1-\beta p}
\left( \int^t_s (t-\xi)^{\alpha-1} \vert v(\xi)\vert d\xi\right)^p
$$
$$
+ C\int^T_0 dt \int^t_0 ds (t-s)^{-1-\beta p}
\left( \int^s_0 \vert (t-\xi)^{\alpha-1} - (s-\xi)^{\alpha-1}\vert 
\vert v(\xi)\vert d\xi\right)^p
=: CI_1 + CI_2.                  \eqno{(2.8)}
$$
\\
\vspace{0.1cm}
\\
{\bf Case 1: $p=1$}.
\\
Equation (2.8) implies 
$$
\int^T_0\int^t_0
\frac{\vert J^{\alpha}v(t) - J^{\alpha}v(s)\vert}
{\vert t-s\vert^{1+\beta}} dsdt
\le CI_1 + CI_2,
$$
where 
$$
I_1 = \int^T_0 \int^t_0 (t-s)^{-1-\beta}
\left( \int^t_s (t-\xi)^{\alpha-1} \vert v(\xi)\vert d\xi\right) dsdt
$$
and
$$
I_2 = \int^T_0 \int^t_0 (t-s)^{-1-\beta}
\left( \int^s_0 \vert (t-\xi)^{\alpha-1} - (s-\xi)^{\alpha-1}\vert
\vert v(\xi)\vert d\xi\right) dsdt.
$$
Exchanging the order of the integral, we obtain
\begin{align*}
& I_1 = \int^T_0 \int^t_0 \vert v(\xi)\vert (t-\xi)^{\alpha-1}
\left( \int^{\xi}_0 (t-s)^{-1-\beta} ds \right) d\xi dt \\
=& \frac{1}{\beta} \int^T_0 \left(
\int^t_0 \vert v(\xi)\vert (t-\xi)^{\alpha-1}
((t-\xi)^{-\beta} - t^{-\beta}) d\xi \right) dt\\
\le &C\int^T_0 \int^t_0 \vert v(\xi)\vert (t-\xi)^{\alpha-\beta-1}
\left( 1 - \left( \frac{t-\xi}{t}\right)^{\beta}\right) d\xi dt
\le C\int^T_0 \left\vert \int^t_0 \vert v(\xi)\vert (t-\xi)^{\ep-1} 
d\xi\right\vert dt
\end{align*}
by $\left\vert \frac{t-\xi}{t} \right\vert \le 1$ and 
$\beta > 0$ for $0<\xi<t$.
Therefore, the Young inequality for the convolution yields
$$
\vert I_1\vert \le C\Vert s^{\xi-1} \, *\, \vert v\vert \Vert_{L^1(0,T)}
\le C\Vert v\Vert_{L^1(0,T)}.
$$
As for the term $I_2$, by the change $s:= \eta(t-\xi) + \xi$ of the variables,
exchanging the order of the integral, we deduce 
\begin{align*}
& \int^t_0 (t-s)^{-1-\beta}\left( \int^s_0
\vert (t-\xi)^{\alpha-1} - (s-\xi)^{\alpha-1}\vert \vert v(\xi)\vert
d\xi \right) ds\\
=& \int^t_0 \left( \int^t_{\xi} (t-s)^{-1-\beta}
\vert (s-\xi)^{\alpha-1} - (t-\xi)^{\alpha-1}\vert ds\right)
\vert v(\xi)\vert d\xi\\
\le &C\int^t_0 \vert v(\xi)\vert (t-\xi)^{\ep-1}
\left( \int^1_0 (1-\eta)^{-1-\beta}(\eta^{\alpha-1} - 1) d\eta \right)d\xi.
\end{align*}
Here for arbitrarily fixed constant $\delta \in \left(0,\, \hhalf\right)$, 
we apply (2.7) to estimate
\begin{align*}
&  \int^1_{1-\delta} (1-\eta)^{-1-\beta}\vert \eta^{\alpha-1} - 1\vert d\eta 
\le \int^1_{1-\delta} (1-\eta)^{-1-\beta}
((1-\alpha)(1-\eta)+C_{\delta}(1-\eta)^2) d\eta\\
\le &C_{\delta} \int^1_{1-\delta} ((1-\eta)^{-\beta} + (1-\eta)^{1-\beta})d\eta
< \infty
\end{align*}
by $0<\beta<\alpha<1$.  Furthermore 
$$
\int^{1-\delta}_0 (1-\eta)^{-1-\beta}\vert \eta^{\alpha-1} - 1\vert d\eta
\le C_{\delta}\int^{1-\delta}_0 (\vert \eta^{\alpha-1}\vert + 1) d\eta
< \infty.
$$
Consequently, the Young inequality yields
$$
\vert I_2\vert \le C_{\delta}\int^T_0 \left\vert \int^t_0 \vert v(\xi)\vert
(t-\xi)^{\ep-1} d\xi \right\vert dt
= C_{\delta}\Vert s^{\ep-1}\, *\, \vert v\vert\Vert_{L^1(0,T)}
\le C_{\delta}\Vert v\Vert_{L^1(0,T)}.
$$
Thus the proof of (i) for $p=1$ and $0<\alpha<1$ is completed.
$\BLS$ 
\\
\vspace{0.1cm}
\\
{\bf Case 2: $p>1$.}
\\
First we note that since $W^{\alpha-\ep,\,p}(0,T) \subset 
W^{\alpha-\ep',\,p}(0,T)$ if $\alpha>\ep' > \ep > 0$, it suffices to 
prove (i) for sufficiently small $\ep>0$, and in particular, we can assume
that 
$$
\alpha > \ep p.                        \eqno{(2.9)}
$$
Let $q \in (1,\infty)$ satisfy $\frac{1}{p} + \frac{1}{q} = 1$.
By $1<p<\infty$, we have $q = \frac{p}{p-1} \in (1,\infty)$.

In $I_1$ and $I_2$, the integrands contain $\vert v(\xi)\vert$, not
$\vert v(\xi)\vert^p$, so that we need to change
such factors to $\vert v(\xi)\vert^p$.  To this end, we resort to 
the H\"older inequality after factorizing 
$$
\vert t-\xi\vert^{\alpha-1} \vert v(\xi)\vert
= (\vert t-\xi\vert^{\alpha-1})^{\frac{1}{q}}
\{(\vert t-\xi\vert^{\alpha-1})^{1-\frac{1}{q}} \vert v(\xi)\vert\}.
$$
Then the H\"older inequality implies
\begin{align*}
& \int^t_s \vert t-\xi\vert^{\alpha-1}\vert v(\xi)\vert d\xi
= \int^t_s ( \vert t-\xi\vert^{\alpha-1})^{\frac{1}{q}}
\{( \vert t-\xi\vert^{\alpha-1})^{1-\frac{1}{q}} \vert v(\xi)\vert\} d\xi\\
\le& \left( \int^t_s (t-\xi)^{\alpha-1} d\xi\right)^{\frac{1}{q}}
\left( \int^t_s ((t-\xi)^{\alpha-1})^{p\left(1-\frac{1}{q}\right)} 
\vert v(\xi)\vert^p d\xi \right)^{\frac{1}{p}},
\end{align*}
and so
\begin{align*}
& \left( \int^t_s \vert t-\xi\vert^{\alpha-1}\vert v(\xi)\vert d\xi
\right)^p
\le \left( \int^t_s (t-\xi)^{\alpha-1} d\xi\right)^{\frac{p}{q}}
\int^t_s \vert t-\xi\vert^{\alpha-1} \vert v(\xi)\vert^p d\xi\\
=& \alpha^{-\frac{p}{q}} (t-s)^{\frac{p}{q}\alpha}
\int^t_s (t-\xi)^{\alpha-1} \vert v(\xi)\vert^p d\xi.
\end{align*}
Therefore, since $\frac{1}{q} = 1 - \frac{1}{p}$ implies
$$
-1-\beta p + \frac{p}{q}\alpha = -1 - \beta p + \alpha p - \alpha
= -1 - \alpha + \ep p,                 \eqno{(2.10)}
$$
using $\int^t_0 \left( \int^t_s \cdots d\xi\right) ds
= \int^t_0 \left( \int^{\xi}_0 \cdots ds \right) d\xi$, 
we obtain
\begin{align*}
& I_1 \le C\int^T_0 dt \int^t_0 (t-s)^{-1-\alpha+\ep p}
\left( \int^t_s (t-\xi)^{\alpha-1}\vert v(\xi)\vert^p d\xi \right)
ds\\
=& C\int^T_0 dt \int^t_0 \left( 
\int^{\xi}_0 (t-s)^{-1-\alpha+\ep p} ds \right) 
(t-\xi)^{\alpha-1}\vert v(\xi)\vert^p d\xi.
\end{align*}
Since $0 < s < \xi < t$, we have
$$
\int^{\xi}_0 (t-s)^{-1-\alpha+\ep p} ds 
= \frac{1}{\ep p -\alpha}((t-\xi)^{-\alpha+\ep p} 
- t^{-\alpha+\ep p}),
$$
and by (2.9) we deduce  
$$
\vert I_1\vert \le C\int^T_0 
\left( \int^t_0 
\vert t^{-\alpha+\ep p} - (t-\xi)^{-\alpha+\ep p}\vert 
\vert t-\xi\vert^{\alpha-1} \vert v(\xi)\vert^p d\xi\right) dt.
                                               \eqno{(2.11)}
$$
On the other hand, 
$$
\vert t^{-\alpha+\ep p} - (t-\xi)^{-\alpha+\ep p}\vert 
= \vert t-\xi\vert^{-\alpha+\ep p}
\left\vert 1 - \left\vert \frac{t-\xi}{t} \right\vert^{\alpha-\ep p}
\right\vert.
$$
From (2.9) and $\left\vert \frac{t-\xi}{t}\right\vert \le 1$, it 
follows that $\left\vert \frac{t-\xi}{t}\right\vert^{\alpha-\ep p} \le 1$
and 
$$
\vert t^{-\alpha+\ep p} - (t-\xi)^{-\alpha+\ep p}\vert
\le 2\vert t-\xi\vert^{-\alpha+\ep p} \quad \mbox{for $0<\xi<t$}.
$$
Hence,
\begin{align*}
& \vert I_1\vert 
\le 2C\int^T_0 \left( \int^t_0 
\vert t-\xi\vert^{-\alpha+\ep p} \vert t-\xi\vert^{\alpha-1}
\vert v(\xi)\vert^p d\xi \right) dt \\
=& C\int^T_0 \left\vert \int^t_0 
\vert t-\xi\vert^{\ep p-1} \vert v(\xi)\vert^p d\xi \right\vert dt
=C\Vert s^{\ep p-1}\,*\, \vert v\vert^p\Vert_{L^1(0,T)}.
\end{align*}
The Young inequality implies 
$$
\vert I_1\vert \le C\Vert s^{\ep p-1}\Vert_{L^1(0,T)}
\Vert \vert v\vert^p\Vert_{L^1(0,T)}
\le C\Vert v\Vert^p_{\LLLP}.                      \eqno{(2.12)}
$$

We now proceed to the proof of $I_2$.  Similarly to $I_1$, we
factorize
\begin{align*}
& \vert (t-\xi)^{\alpha-1} - (s-\xi)^{\alpha-1}\vert \vert v(\xi)\vert\\
=& \vert (t-\xi)^{\alpha-1} - (s-\xi)^{\alpha-1}\vert^{\frac{1}{q}}
\{ \vert (t-\xi)^{\alpha-1} - (s-\xi)^{\alpha-1}\vert^{1-\frac{1}{q}}
\vert v(\xi)\vert\}
\end{align*}
and apply the H\"older inequality to obtain
\begin{align*}
& \left( \int^s_0 \vert (t-\xi)^{\alpha-1} - (s-\xi)^{\alpha-1}\vert 
\vert v(\xi)\vert d\xi\right)^p\\
\le & \left( \int^s_0 \vert (t-\xi)^{\alpha-1} - (s-\xi)^{\alpha-1}\vert 
d\xi\right)^{\frac{p}{q}}
\int^s_0 \vert (t-\xi)^{\alpha-1} - (s-\xi)^{\alpha-1}\vert 
\vert v(\xi)\vert^p d\xi\\
=& \left( \frac{(t-s)^{\alpha} - (t^{\alpha}-s^{\alpha})}{\alpha}
\right)^{\frac{p}{q}}
\int^s_0 \vert (t-\xi)^{\alpha-1} - (s-\xi)^{\alpha-1}\vert 
\vert v(\xi)\vert^p d\xi\\
\le& C(t-s)^{\frac{p\alpha}{q}}
\int^s_0 \vert (t-\xi)^{\alpha-1} - (s-\xi)^{\alpha-1}\vert 
\vert v(\xi)\vert^p d\xi.
\end{align*}
Therefore, 
\begin{align*}
& \int^t_0 (t-s)^{-1-\beta p}
\left( \int^s_0 \vert (t-\xi)^{\alpha-1} - (s-\xi)^{\alpha-1}\vert 
\vert v(\xi)\vert d\xi \right)^p ds\\
\le& C\int^t_0 (t-s)^{-1-\beta p+\frac{p\alpha}{q}}
\left( \int^s_0 \vert (t-\xi)^{\alpha-1} - (s-\xi)^{\alpha-1}\vert 
\vert v(\xi)\vert^p d\xi \right) ds.
\end{align*}
Exchanging the order of the integral and using (2.10), we have
\begin{align*}
& \int^t_0 (t-s)^{-1-\beta p}
\left( \int^s_0 \vert (t-\xi)^{\alpha-1} - (s-\xi)^{\alpha-1}\vert 
\vert v(\xi)\vert^p d\xi \right) ds\\
\le& C\int^t_0 \left( \int^t_{\xi} (t-s)^{-1-\alpha + p\ep}
((s-\xi)^{\alpha-1} - (t-\xi)^{\alpha-1}) ds \right) 
\vert v(\xi)\vert^p d\xi.
\end{align*}
Changing the variables $s \mapsto \eta$ by $s:= \eta(t-\xi)+\xi$, we obtain
\begin{align*}
& \int^t_0 (t-s)^{-1-\alpha+\ep p}
\left( \int^s_0 \vert (t-\xi)^{\alpha-1} - (s-\xi)^{\alpha-1}\vert 
\vert v(\xi)\vert^p d\xi \right) ds\\
= &C\int^t_0 \vert v(\xi)\vert^p (t-\xi)^{\ep p -1}
\left( \int^1_0 (1-\eta)^{-1-\alpha+\ep p}
(\eta^{\alpha-1}-1) d\eta \right) d\xi.
\end{align*}
Here 
$$
\int^1_0 (1-\eta)^{-1-\alpha+\ep p}(\eta^{\alpha-1}-1) d\eta < \infty.
$$
Indeed, for fixed $\delta \in \left(0, \hhalf\right)$, we deduce 
$$
\int^{\delta}_0 (1-\eta)^{-1-\alpha+\ep p}(\eta^{\alpha-1}-1) d\eta
\le \int^{\delta}_0 (\eta^{\alpha-1}-1) d\eta < \infty
$$
and (2.7) yields
\begin{align*}
& \int^1_{1-\delta} (1-\eta)^{-1-\alpha+\ep p}
\vert \eta^{\alpha-1} - 1\vert d\eta \\
\le & \int^1_{1-\delta} (1-\eta)^{-1-\alpha+\ep p}
((1-\alpha)(1-\eta) + C_{\delta}(1-\eta)^2) d\eta
\le C_{\delta}\int^1_{1-\delta} (1-\eta)^{-\alpha+\ep p} d\eta < \infty,
\end{align*}
because $-\alpha + \ep p > -1$ by $0<\alpha<1$.

Therefore, 
\begin{align*}
& \int^t_0 (t-s)^{-1-\alpha + \ep p}
\left( \int^s_0 \vert (t-\xi)^{\alpha-1} - (s-\xi)^{\alpha-1}\vert 
\vert v(\xi)\vert^p d\xi \right) ds\\
\le& C\vert (s^{\ep p-1}\, *\, \vert v\vert^p)(t)\vert, \quad 0<t<T.
\end{align*}
Then the Young inequality yields
$$
\vert I_2\vert \le C\Vert s^{\ep p-1}\, *\, \vert v\vert^p\Vert
_{L^1(0,T)} \le C\Vert s^{\ep p -1}\Vert_{L^1(0,T)}
\Vert \vert v\vert^p\Vert_{L^1(0,T)} \le C\Vert v\Vert_{\LLLP}^p.
                                     \eqno{(2.13)}
$$
Hence, for $p > 1$ and $0<\alpha<1$, by (2.8), (2.12) and (2.13), 
we complete the proof of the estimate in Theorem 2.1 (i).
\\
\vspace{0.1cm}
\\
{\bf Second Step: $\alpha > 1$.}
\\
Let $\alpha = \ell + \sigma$ where $\ell \in \N$ and $0<\sigma<1$.
Then
$$
J^{\alpha}v(t) = \frac{1}{\Gamma(\ell+\sigma)}\int^t_0
(t-s)^{\ell+\sigma-1} v(s) ds, \quad 0<t<T, \, v\in \LLLP.
$$
Therefore 
$$
\frac{d^k}{dt^k}J^{\alpha}v(t) = \frac{1}{\Gamma(\ell-k+\sigma)}\int^t_0
(t-s)^{\ell-k+\sigma-1} v(s) ds, \quad 0\le k\le \ell, \,
0<t<T, \, v\in \LLLP.
$$
In particular, 
$$
\frac{d^{\ell}}{dt^{\ell}}J^{\alpha}v(t) = \frac{1}{\Gamma(\sigma)}\int^t_0
(t-s)^{\sigma-1} v(s) ds, \quad 0\le k\le \ell, \,
0<t<T, \, v\in \LLLP.
$$
Hence,
$$
\left\Vert \frac{d^k}{dt^k}J^{\alpha}\right\Vert_{\LLLP}
= \frac{1}{\Gamma(\ell-k+\sigma)}\Vert s^{\ell-k+\sigma-1}\,*\,
v\Vert_{\LLLP} \le C\Vert v\Vert_{\LLLP}, \quad 0\le k\le \ell.
                                                      \eqno{(2.14)}
$$
Since we have already proved (2.5) with $\alpha \in (0,1)$, 
we see
$$
\left\Vert \frac{1}{\Gamma(\sigma)}\int^t_0
(t-s)^{\sigma-1} v(s) ds\right\Vert_{W^{\sigma-\ep,p}(0,T)}
\le C\Vert v\Vert_{\LLLP}                       \eqno{(2.15)}
$$
with small $\ep \in (0, \sigma)$.   
Consequently, in terms of (2.14) and (2.15), we deduce
\begin{align*}
& \Vert J^{\alpha}v\Vert^p_{W^{\alpha-\ep,p}(0,T)}
= \sum_{k=0}^{\ell} \left\Vert \frac{d^k}{dt^k}J^{\alpha}v\right\Vert^p
_{\LLLP}\\
+& \int^T_0\int^T_0 \frac{
\left\vert \frac{d^{\ell}}{dt^{\ell}}J^{\alpha}v(t)
- \frac{d^{\ell}}{dt^{\ell}}J^{\alpha}v(s)\right\vert^p}
{\vert t-s\vert^{1+(\sigma-\ep)p}}dtds
\le C\Vert v\Vert^p_{\LLLP}
\end{align*}
for all $v\in \LLLP$.

Finally, for $\alpha = m-1+\sigma$ with $m\in \N$ and
$0<\sigma<1$, we have to prove that 
$\WWW \subset \, _{0}W^{\alpha-\ep,p}(0,T)$.
Let $u \in \WWW$ be arbitrarily given.  By the definition (2.4), there exists
$w \in \LLLP$ such that $u=J^{\alpha}w$.  By the density, we can find
$w_n\in \, _{0}C^m[0,T]$, $n\in \N$, such that $w_n \RRRR w$ in 
$\LLLP$ as $n \to \infty$.  Setting $u_n:= J^{\alpha}w_n$, $n\in \N$,
by $w_n \in \, _{0}C^m[0,T]$, we can directly prove that 
$u_n \in \, _{0}C^m[0,T]$ for $n\in \N$.
\\
Indeed, since we can directly see that $J^{\alpha}C[0,T] \subset C[0,T]$
with $\alpha>0$ by estimating $J^{\alpha}u(t) - J^{\alpha}u(s)$ for 
$u \in C[0,T]$, this follows by applying (2.16) to $v:= w_n$ and 
$\gamma:=\alpha$:  
$$
\frac{d^k}{dt^k}J^{\gamma}v = J^{\gamma}\frac{d^kv}{dt^k}
\quad \mbox{for $v \in \CCCM$ and $\gamma>0$, $k=0,1,..., m$.}
                                             \eqno{(2.16)}
$$
\\
{\bf Verification of (2.16).}
\\
Since $\frac{d^kv}{dt^k}(0) = 0$ for $0\le k \le m-1$
and $\int^t_0 (t-s)^{\gamma-1}v(s) ds 
= \int^t_0 s^{\gamma-1}v(t-s) ds$, we have
$$
\frac{d}{dt}J^{\gamma}v(t)
= \frac{1}{\Gamma(\gamma)}\left( \int^t_0 s^{\gamma-1}\frac{dv}{dt}(t-s) ds 
+ t^{\gamma-1}v(0) \right)
= \frac{1}{\Gamma(\gamma)}\int^t_0 s^{\gamma-1}\frac{dv}{dt}(t-s) ds 
$$
and
$$
\frac{d^2}{dt^2}J^{\gamma}v(t)
= \frac{1}{\Gamma(\gamma)}\left( \int^t_0 s^{\gamma-1}\frac{d^2v}{dt^2}(t-s) ds + t^{\gamma-1}\frac{dv}{dt}(0) \right).
$$
We can continue the calculations to finish the proof of (2.16).
$\BLS$
\\
Then, by (2.5), we can deduce that  
$J^{\alpha}w_n \RRRR J^{\alpha}w$ in $W^{\alpha-\ep,p}(0,T)$, that is,
$u_n \RRRR u$ in $W^{\alpha-\ep,p}(0,T)$ as $n\to \infty$.
Since $u_n \in \, _{0}C^m[0,T]$, this means that 
$u \in \, \ooo{_{0}C^m[0,T]}^{W^{\alpha-\ep,p}(0,T)}
= \, _{0}W^{\alpha-\ep,p}(0,T)$.  

Thus the proof of Theorem 2.1 (i) is completed.
$\BLS$
\\
\vspace{0.1cm}
\\
{\bf Proof of Theorem 2.1 (ii).}
\\
{\bf First Step.}
\\
We show
\\
{\bf Lemma 2.4.}
\\
{\it
Let $\alpha = m-1+\sigma$ with $m\in \N$ and $0<\ep<1-\sigma$, and
$1 \le p < \infty$.  Then there exists a constant $C=C(p,\alpha,\ep)>0$ such 
that 
$$
\Vert J^{m-\alpha}u\Vert_{W^{m,p}(0,T)} \le C\Vert u\Vert
_{W^{\alpha+\ep,p}(0,T)}
$$
for all $u \in \,_{0}W^{\alpha+\ep,p}(0,T)$.
}

For $p=1$ and $m=1$, the proof is found in Theorem 4.2.3 (pp.77-78) in 
\cite{GV}.  Also our proof for the general $m\in \N$ and $1<p<\infty$,
relies on direct estimation, but is lengthy.
Thus the proof of Lemma 2.4 is postponed to Section 6.
\\
{\bf Second Step.}
\\
By (2.16), we can directly prove that formula (2.6) holds for 
$u \in \CCCM$. 
Indeed, by $m-\alpha = 1-\sigma \in (0,1)$, equality (2.16) yields
$$
\frac{d^m}{dt^m}J^{m-\alpha}u = J^{m-\alpha}\frac{d^mu}{dt^m},
$$
and so 
$$
J^{\alpha}\left( \frac{d^m}{dt^m}J^{m-\alpha}u\right)
= J^{\alpha}J^{m-\alpha}\frac{d^mu}{dt^m}
= \frac{1}{\Gamma(m)}\int^t_0 (t-s)^{m-1} \frac{d^mu}{ds^m}(s) ds.
$$
Here we used $J^{\alpha}J^{m-\alpha} = J^m$ by Lemma 2.1 (i).
Repeating the integration by parts and using 
$\frac{d^ku}{ds^k}(0) = 0$ for $0\le k \le m-1$. we see 
$$
\frac{1}{\Gamma(m)}\int^t_0 (t-s)^{m-1} \frac{d^mu}{ds^m}(s) ds
= u(t), \quad 0<t<T.
$$
Thus we have verified (2.6) for $u \in \CCCM$.
\\

Next let $u \in \, _{0}W^{\alpha+\ep,p}(0,T)$ be arbitrarily given.
By the definition, we can choose $u_n \in \CCCM$, $n\in \N$, such that 
$u_n \RRRR u$ in $W^{\alpha+\ep,p}(0,T)$ as $n\to \infty$.
Applying Lemma 2.4, we see that 
$J^{m-\alpha}u_n \RRRR J^{m-\alpha}u$ in $W^{m,p}(0,T)$, that is,
$\frac{d^m}{dt^m}J^{m-\alpha}u_n \, \RRRR\, \frac{d^m}{dt^m}J^{m-\alpha}u$
in $\LLLP$ as $n\to \infty$ and $\frac{d^m}{dt^m}J^{m-\alpha}u
\in \LLLP$.  Therefore, Lemma 2.2 yields
$$
J^{\alpha}\left( \frac{d^m}{dt^m}J^{m-\alpha}u_n\right)
\, \RRRR\, J^{\alpha}\left( \frac{d^m}{dt^m}J^{m-\alpha}u\right)
\quad \mbox{in $\WWW$}
$$
as $n\to \infty$, that is, by (2.6) with $u_n \in \, _{0}C^m[0,T]$,
we obtain
$$
u_n \, \RRRR\, J^{\alpha}\left( \frac{d^m}{dt^m}J^{m-\alpha}u\right)
\quad \mbox{in $\WWW$}
$$
as $n\to \infty$.  
Consequently, since $u_n \RRRR u$ in $W^{\alpha+\ep,p}(0,T)$
as $n\to \infty$, we verify (2.6) for all 
$u\in \, _{0}W^{\alpha+\ep,p}(0,T)$.
Thus the proof of Theorem 2.1 (ii) is complete.
$\BLS$
\\
\vspace{0.1cm}
\\
{\bf Proof of Theorem 2.1 (iii).}
\\
Since the embedding $W^{\alpha+\ep,p}(0,T)\,\RRRR \, \LLLP$ is compact
(e.g., \cite{Ad}, \cite{Gri}), part (i) implies
$$
\mbox{$J^{\alpha}: \LLLP \, \RRRR\, \LLLP$ is a compact operator.}
                                               \eqno{(2.17)}
$$
This compactness is known (e.g., \cite{GV}) but it follows from (i).

Let $u_n \in W_{\alpha+\ep,p}(0,T)$, $n\in \N$ satisfy 
$$
\sup_{n\in \N} \Vert u_n\Vert_{W_{\alpha+\ep,p}(0,T)} < \infty.
                                              \eqno{(2.18)}
$$
Then, for each $n\in \N$, there exists $w_n\in \LLLP$ such that 
$u_n=J^{\alpha+\ep}w_n$.  Moreover the definition of the norm implies
$\Vert u_n\Vert_{W_{\alpha+\ep,p}(0,T)} = \Vert w_n\Vert_{\LLLP}$ for 
$n\in \N$.  Hence,
$$
\sup_{n\in \N} \Vert w_n\Vert_{\LLLP} < \infty.
                                              \eqno{(2.19)}
$$
Setting $v_n := J^{\ep}w_n$, the Young inequality yields
$$
\Vert v_n\Vert_{\LLLP} \le C\Vert w_n\Vert_{\LLLP}.
$$
Therefore, in view of (2.17) and (2.19), we conclude that $\{ v_n\}_{n\in\N}$
contains a convergent subsequence $\{ v_{n(k)}\}_{k\in \N}$ in 
$\LLLP$.  Since $u_{n(k)} = J^{\alpha}J^{\ep}w_{n(k)}
= J^{\alpha}v_{n(k)}$, by the definition of 
the norm $\Vert \cdot\Vert_{\WWW}$. we see that $u_{n(k)}$ is convergent 
in $\WWW$.  Thus by (2.18), we complete the proof of part (iii) and thus
the proof of Theorem 2.1 is complete.
$\BLS$
\\
\vspace{0.1cm}
\\
{\bf Proof of Corollary 2.1.}
\\
(i) We assume that 
$$
\sup_{n\in \N} \Vert u_n\Vert_{W_{\beta,p}(0,T)} < \infty.
                                              \eqno{(2.20)}
$$
By the definition of $W_{\beta,p}(0,T)$. we can find $w_n \in \LLLP$
such that $u_n=J^{\beta}w_n$ for $n\in \N$.
Moreover $\Vert w_n\Vert_{\LLLP} = \Vert u_n\Vert_{W_{\beta,p}(0,T)}$
for $n\in \N$, so that (2.20) implies 
$$
\sup_{n\in \N} \Vert w_n\Vert_{\LLLP} < \infty.
                                              \eqno{(2.21)}
$$
Since $J^{\alpha}u_n = J^{\alpha}J^{\beta}w_n$, by the definition of 
$\Vert \cdot\Vert_{W_{\beta,p}(0,T)}$, we deduce 
$$
\Vert J^{\alpha}u_n\Vert_{W_{\beta,p}(0,T)} 
= \Vert J^{\beta}(J^{\alpha}w_n)\Vert_{W_{\beta,p}(0,T)}
= \Vert J^{\alpha}w_n\Vert_{\LLLP}.
$$
In terms of (2.17) and (2.21), there exists a sequence $\{n(k)\}
_{k\in \N} \subset \N$ such that 
$J^{\alpha}w_{n(k)}$ is convergent in $\LLLP$, that is,
$J^{\alpha}u_{n(k)}$ is convergent in $W_{\beta,p}(0,T)$.
With (2.21), we deduce that $J^{\alpha}: W_{\beta,p}(0,T)\, 
\RRRR\, W_{\beta,p}(0,T)$ is a compact operator.
$\BLS$
\\
(ii) If $u \in W_{\beta,p}(0,T)$, then $u=J^{\beta}w$ with 
$w\in \LLLP$.  By Lemma 2.1 (i), we see that $J^{\alpha}u =  
J^{\alpha}J^{\beta}w = J^{\alpha+\beta}w\in 
W_{\alpha+\beta,p}(0,T)$.  
Therefore, $J^{\alpha}W_{\beta,p}(0,T) \subset 
W_{\alpha+\beta,p}(0,T)$.

Conversely, let $u \in W_{\alpha+\beta,p}(0,T)$.  Then
there exists $w\in \LLLP$ such that $u = J^{\alpha+\beta}w
= J^{\alpha}(J^{\beta}w)$.  By $J^{\beta}w \in W_{\beta,p}(0,T)$, we 
see that $u \in J^{\alpha}W_{\beta,p}(0,T)$, that is,
$W_{\alpha+\beta,p}(0,T) \subset J^{\alpha}W_{\beta,p}(0,T)$.

Finally we will prove the norm equivalence.  Let $u\in W_{\beta,p}(0,T)$
be arbitrarily given.
Then $u=J^{\beta}w$ with some $w\in \LLLP$ and 
$\Vert u\Vert_{W_{\beta,p}(0,T)} = \Vert w\Vert_{\LLLP}$.
Moreover the definition of $\Vert \cdot\Vert_{W_{\alpha+\beta,p}(0,T)}$
implies
$$
\Vert J^{\alpha}u\Vert_{W_{\alpha+\beta,p}(0,T)}
= \Vert J^{\alpha}J^{\beta}w\Vert_{W_{\alpha+\beta,p}(0,T)}
= \Vert J^{\alpha+\beta}w\Vert_{W_{\alpha+\beta,p}(0,T)}
= \Vert w\Vert_{\LLLP}.
$$
and so $\Vert J^{\alpha}u\Vert_{W_{\alpha+\beta,p}(0,T)}
= \Vert u\Vert_{W_{\beta,p}(0,T)}$.
Thus the proof of Corollary 2.1 is complete.
$\BLS$
\\

We close this section with the following proposition which is used
in Section 4 and is proved easily.
\\
{\bf Proposition 2.1.}
\\
{\it
(i) $W_{\alpha,p}(0,T) \subset W_{\beta,p}(0,T)$ if $\alpha > \beta > 0$ and 
there exists a constant $C>0$ such that 
$$
\Vert u\Vert_{W_{\beta,p}(0,T)} \le C\Vert u\Vert_{W_{\alpha,p}(0,T)}
$$
for each $u \in W_{\alpha,p}(0,T)$.
\\
(ii) $J^{\alpha}\ppp_t^{\beta}u = J^{\alpha-\beta}u$ for $\alpha > \beta > 0$
and each $u \in W_{\beta,p}(0,T)$.
}
\\
{\bf Proof of Proposition 2.1.}
\\
(i) Let $u \in W_{\alpha,p}(0,T)$.  Then $u = J^{\alpha}w$ with some 
$w\in \LLLP$.  Since $J^{\alpha}w = J^{\beta}J^{\alpha-\beta}w$ 
by Lemma 2.1 (i), we have $u = J^{\alpha}w = J^{\beta}(J^{\alpha-\beta}w)
\in J^{\beta}\LLLP$.
Moreover, by $u = J^{\alpha}w = J^{\beta}(J^{\alpha-\beta}w)$, 
the definition of
the norm implies that 
$\Vert u\Vert_{W_{\alpha,p}(0,T)} = \Vert w\Vert_{\LLLP}$ and 
$$
\Vert u\Vert_{W_{\beta,p}(0,T)} = \Vert J^{\alpha-\beta}w\Vert_{\LLLP}
\le C\Vert u\Vert_{\LLLP}.
$$
Therefore, $\Vert u\Vert_{W_{\beta,p}(0,T)} \le C\Vert u\Vert
_{W_{\alpha,p}(0,T)}$.   Thus part (i) is proved.
$\BLS$
\\
(ii) Let $u \in W_{\beta,p}(0,T)$.  Then $u = J^{\beta}w$ with some
$w\in \LLLP$ and
$$
J^{\alpha}\ppp_t^{\beta} = J^{\alpha}(J^{\beta})^{-1}J^{\beta}w
= J^{\alpha}w.
$$
On the other hand, by $\alpha - \beta, \beta > 0$, Lemma 2.1 (i) implies
$$
J^{\alpha-\beta}u = J^{\alpha-\beta}J^{\beta}w 
= J^{(\alpha-\beta)+\beta}w = J^{\alpha}w.
$$
Hence $J^{\alpha}\ppp_t^{\beta} = J^{\alpha-\beta}u$ for $u \in 
W_{\beta,p}(0,T)$.
Thus the proof of Proposition 2.1 is complete.
$\BLS$
\section{Fractional derivative $\pppa$ in $L^p(0,T)$}

We consider an initial value problem (1.2) and (1.3):
$$
\ddda u(t) = b(t)u(t) + f(t), \quad 0<t<T, \qquad u(0) = a,
$$
where $f$ is singular in the sense that $f \not\in L^r(0,T)$ for any $r \ge 1$.
For example, let $f(t)$ be the Dirac delta function which means an 
impulsive source at $t=0$.  It is desirable to construct a framework
in order to treat such singular terms in fractional differential equations.
For the mathematical treatments, we need the 
formulation of time-fractional derivative $\pppa$ in Sobolev spaces of 
non-positive orders.
In Yamamoto \cite{Ya21}, we find such studies in the space $L^2(0,T)$ and 
the treatment with fixed $p=2$ is less flexible.

In this section, we define $\pppa$ in $L^p(0,T)$-based Sobolev-Slobodecki 
spaces of negative orders.  Our extension of the domain of $\pppa$ 
from $\WWW$, relies on 
the adjoint of fractional differential operator.  
\\

We set
$$
J_{\alpha}v(t) = \frac{1}{\Gamma(\alpha)}\int^T_t (\xi-t)^{\alpha-1} v(\xi) 
d\xi, \quad 0<t<T, \, v \in L^1(0,T).        \eqno{(3.1)}
$$
Then, similarly to Lemma 2.1, we can prove
\\
{\bf Lemma 3.1.}
\\
{\it
Let $\alpha, \beta > 0$.
\\ 
(i) $J_{\alpha}J_{\beta}v = J_{\alpha+\beta}v$ for 
$v \in L^1(0,T)$.
\\
(ii) $J_{\alpha}: \LLLP \, \RRRR\, \LLLP$ is an injective compact operator.
}

We can prove directly but we can transform the corresponding results for 
$J^{\alpha}$ through the following transformation:
$$
\tau:\, \LLLP \, \RRRR\, \LLLP, \qquad 
(\tau v)(t) = v(T-t), \quad 0<t<T.            \eqno{(3.2)}
$$
For example, as is readily verified,
$$
J_{\alpha}v(t) = (\tau J^{\alpha}(\tau v))(t), \quad 0<t<T,\,
v\in \LLLP.                \eqno{(3.3)}
$$
In view of (3.3), Lemma 3.1 follows directly from Lemma 2.1.

In the same way as in Section 2, for $1\le p < \infty$ and $\alpha>0$, by 
Lemma 3.1, we can define
$$
\WWWR := J_{\alpha}\LLLP, \quad
\Vert u\Vert_{\WWWR}:= \Vert (J_{\alpha})^{-1}u\Vert_{\LLLP} \quad 
\mbox{for $u \in \, _{\alpha,p}W(0,T)$}.           \eqno{(3.4)}
$$

Moreover, we can similarly prove
\\
{\bf Lemma 3.2.}
\\
{\it
Let $\alpha > 0$.  If 
$w_n \RRRR w$ in $\LLLP$, then $J_{\alpha}w_n \RRRR J_{\alpha}w$ in 
$\WWWR$.
}
\\
{\bf Lemma 3.3.}
\\
{\it
The space $\WWWR$ is a Banach space with the norm 
$\Vert \cdot\Vert_{\WWWR}$.
}
\\

For $m\in \N$, we set  
$$
^{0}C^m[0,T]:= 
\left\{ v\in C^m[0,T];\, \frac{d^kv}{dt^k}(T) = 0\quad 
\mbox{for $k=0,1,..., m-1$}\right\}
$$
and 
$$
\WWWRO := \ooo{\, ^{0}C^m[0,T]}^{W^{\alpha,p}(0,T)}.
$$
We remark that $\,^{0}C^m[0,T] = \tau \CCCM$ and
$\tau: \, \WWW\, \RRRR \, \WWWR$ is an isomorphism.

Then, by means of Theorem 2.1 and the symmetric transform $\tau$
in $t$, we can readily prove
\\
{\bf Theorem 3.1.}
\\
{\it
Let $\alpha > 0$, $\not\in \N$ and $0<\ep<\alpha$.
\\
(i) $\, ^{0}W^{\alpha+\ep,p}(0,T) \, \subset \, \WWWR\,
\subset \, ^{0}W^{\alpha-\ep,p}(0,T)$.
Moreover there exists a constant $C=C(\alpha,\ep)>0$ such that 
$$
\Vert J_{\alpha}v\Vert_{W^{\alpha-\ep,p}(0,T)}
\le C\Vert v\Vert_{\LLLP} \quad \mbox{for all $v\in \LLLP$.}
$$
\\
(ii) The embedding $\, _{\alpha+\ep,p}W(0,T)\, \RRRR\, \WWWR$ is 
compact.
\\
(iii) For $\alpha>0$ and $\beta \ge 0$, the operator
$J_{\alpha}:\, _{\beta,p}W(0,T)\, \RRRR\,  _{\beta,p}W(0,T)$ is 
compact.
}
\\

We will use $J_{\alpha}$ in order to define $\pppa$ in Sobolev spaces of 
negative orders through he adjoint.  As for notations and terminologies,
we follow Brezis \cite{Bre}.

Let $\WWWR^*$ denote the dual space of $\WWWR$, that is, the space of all the 
bounded linear real-valued functionals defined over $\WWWR$.

In place of $u(\psi)$, by $\,_{\WWWR^*}<u,\, \psi>_{\WWWR}$, we denote 
the value of $u$ at
$\psi \in \WWWR$.  As the norm $\Vert \cdot\Vert_{\WWWR^*}$ we follow a usual 
definition:
$$
\Vert u\Vert_{\WWWR^*}\, := \, \sup_{\Vert \psi\Vert_{\WWWR}\,=1}
\vert \,_{\WWWR^*}<u,\, \psi>_{\WWWR}\vert.
$$

Henceforth, throughout this article, for $1\le p < \infty$,
we define $1 < q \le \infty$ by 
$$
q = 
\left\{ \begin{array}{rl}
\frac{p}{p-1} \quad &\mbox{if $p>1$}, \\
\infty, \quad &\mbox{if $p=1$},
\end{array}\right.
$$
and we set 
$$
_{L^q(0,T)}(v, \, u)_{\LLLP} := \int^T_0 u(t)v(t) dt, \quad
v\in L^q(0,T), \, u\in \LLLP \quad \mbox{with $1\le p < \infty$}.
$$
We note that $\frac{1}{p} + \frac{1}{q} = 1$ and 
$$
L^p(0,T)^* = L^q(0,T) \quad \mbox{if $1\le p < \infty$}
$$
(e.g., \cite{Bre}).
\\
{\bf Example of $\WWWR^*$.}
\\  
Let $\alpha p > 1$.  Choosing $\ep>0$ sufficiently small such 
that $(\alpha-\ep)p>1$, by the Sobolev embedding 
(e.g., \cite{Ad}, \cite{Gri}) by Theorem 3.1 (i) we verify that 
$\WWWR \subset C[0,T]$.
Therefore, for $t_0\in [0,T]$, we define $\delta_{t_0} \in \WWWR^*$ by 
$$
\, _{\WWWR^*}<\delta_{t_0},\, \psi>_{\WWWR}\, : = \psi(t_0)
$$
and we can see $\delta_{t_0} \in \WWWR^*$ if $\alpha p > 1$.
\\

We here sum up useful results.
\\
{\bf Lemma 3.4.}
\\
{\it
(i) $C^{-1}\Vert J_{\alpha}u\Vert_{\LLLP} 
\le \Vert u\Vert_{\LLLP} \le C\Vert J_{\alpha}u\Vert
_{\WWWR}$ for each $u \in \LLLP$.
\\
(ii) For arbitrarily chosen $v \in \LLLQ$, we define a mapping:
$$
F_v:\, u \mapsto \, _{\LLLQ}(v, \, u)_{\LLLP} \quad \mbox{for each 
$u \in \, _{\alpha,p}W(0,T)$}. 
$$
Then, $F_v$ is a bounded linear functional on $\, _{\alpha,p}W(0,T)$.
Identifying $F_v \in \, _{\alpha,p}W(0,T)^*$ with $v$, in place of
$F_v(u)$ we write 
$$
_{\, _{\alpha,p}W(0,T)^*}<v, \, u>_{\, _{\alpha,p}W(0,T)}
\,\, = \, _{\LLLQ}(v,\, u)_{\LLLP},
$$
so that 
$$
\Vert F_v\Vert_{\WWWR^*} = \Vert v\Vert_{\WWWR^*}
\le C\Vert v\Vert_{\LLLQ}.
$$
\\
(iii) 
$$
J_{\alpha}^*:\, _{\alpha+\beta,p}W(0,T)^* \, \RRRR\, 
_{\beta,p}W(0,T)^*
$$
is an isomorphism.  In particular, 
$$
\Vert J_{\alpha}^*u\Vert_{_{\beta,p}W(0,T)^*} 
= \Vert u\Vert_{_{\alpha+\beta}W(0,T)^*} \quad 
\mbox{for all $u \in \, _{\alpha+\beta}W(0,T)^*$.}
$$
}
\\
{\bf Proof of Lemma 3.4.}
\\
(i) By the defintion of the norm $\Vert \cdot\Vert_{\WWWR}$,  
we have $\Vert J_{\alpha}u\Vert_{\WWWR} = \Vert u\Vert_{\LLLP}$
for any $u \in \LLLP$.
The Young inequality implies $C\Vert J_{\alpha}u\Vert_{\LLLP}
\le \Vert u\Vert_{\LLLP}$.  Thus the proof of the lemma (i) is 
complete.
$\BLS$
\\
(ii) Immediately we see that $F_v$ is well-defined for 
$u \in \WWWR$ and is a lineat mapping.
Setting $u:= J_{\alpha}w$ with $w\in\LLLP$, we have 
$\Vert u\Vert_{\WWWR} = \Vert J_{\alpha}w\Vert_{\WWWR}
= \Vert w\Vert_{\LLLP}$, by Lemma 3.4 (i) we obtain
\begin{align*}
& \Vert v\Vert_{\WWWR^*}
= \, \sup_{\Vert u\Vert_{\WWWR}\, = 1} \vert 
\, _{\LLLQ}(v, \, u)_{\LLLP}\vert 
= \, \sup_{\Vert w\Vert_{\LLLP}\, = 1} \vert 
\, _{\LLLQ}(v,\, J_{\alpha}w)_{\LLLP}\vert\\
\le& \Vert v\Vert_{\LLLQ} \sup_{\Vert w\Vert_{\LLLP}\, = 1} 
\Vert J_{\alpha}w\Vert_{\LLLP} 
\le C\Vert v\Vert_{\LLLQ}.
\end{align*}
Hence, $\Vert v\Vert_{\WWWR^*} \le C\Vert v\Vert_{\LLLQ}$.
Thus the proof of (ii) is complete.
$\BLS$
\\
(iii) Thanks to the operator $\tau$ defined by (3.2), 
Corollary 2.1 (ii) yields that
$$
J_{\alpha}:\, _{\beta,p}W(0,T)\, \RRRR \, _{\alpha+\beta,p}W(0,T)
$$
is an isomorphism for $\alpha > 0$ and $\beta \ge 0$.
Consequently, part (iii) follows directly by the closed range theorem 
(e.g., Section 7 in Chapter 2 of \cite{Bre}).
Thus the proof of Lemma 3.4 is complete.
$\BLS$
\\

For arbitrarily chosen $v\in \LLLQ$ and $u \in \LLLP$, we apply 
Lemma 3.4 (ii) by setting $u:= J_{\alpha}u$, we see that
$$
_{\WWWR^*}<v, \, J_{\alpha}u>_{\WWWR}\, := \, _{\LLLQ}(v, \, J_{\alpha}u)
_{\LLLP}                                                    \eqno{(3.5)}
$$
for all $u\in \LLLP$ and $v\in \LLLQ$.
                          
On the other hand, exchanging the order of the integral
$\int^t_0 \left(\int^s_0 \cdots d\xi\right) ds
= \int^t_0 \left( \int^t_{\xi} \cdots ds \right)d\xi$, we can readily verify  
$$
\, _{\LLLQ}(J^{\alpha}v, \, u)_{\LLLP} 
\, = \, _{\LLLQ}(v, \, J_{\alpha}u)_{\LLLP} \quad \mbox{for all 
$u\in \LLLP$ and $v\in \LLLQ$}.
$$
In terms of (3.5), we can reach 
$$
_{\WWWR^*}<v, \, J_{\alpha}u>_{\WWWR} \,\,
= \, _{\LLLQ}(J^{\alpha}v, \, u)_{\LLLP}           \eqno{(3.6)}
$$
for all $u\in \LLLP$ and $v\in \LLLQ$.
We note that since $J_{\alpha}\,:\LLLP \, \RRRR\, \WWWR$ is 
bounded, the adjoint operator $J_{\alpha}^*\,: 
\WWWR^*\, \RRRR\, L^p(0,T)^* = L^q(0,T)$ exists and 
$J_{\alpha}^*$ is the operator with the maximal domain among 
operators $J: \, \WWWR^*\, \RRRR\, L^q(0,T)$ satisfying 
$$
_{\WWWR^*}<v, \, J_{\alpha}u>_{\WWWR}\, = \, _{\LLLQ}(Jv, \, u)_{\LLLP} 
                                            \eqno{(3.7)}
$$
for all $v \in \DDD(J) \subset \, \WWWR^*$ and $u\in \LLLP$.   
By taking $J:= J^{\alpha}$ with $\DDD(J) = \LLLQ$, from (3.6) it follows that 
(3.7) holds.  Therefore, the maximality of $J_{\alpha}^*$ yields 
$$
J^{\alpha} \subset J_{\alpha}^*, \quad 
\DDD(J_{\alpha}^*) = \, \WWWR\, \supset \LLLQ.       \eqno{(3.8)}
$$

Now we define 
\\
{\bf Definition 3.1.}
\\
{\it We define
$$
\pppa = (J_{\alpha}^*)^{-1}, \quad \DDD(\pppa) = \,\WWWR^*.
$$
}
\\

In view of (3.8), we extend $\pppa$ defined on $\WWW$ to the domain 
$\, \WWWR^*$, which means that $\pppa$ in Definitions 2.1 and 3.1 
coincide in $\WWW$: 
$$
\pppa u = (J^{\alpha})^{-1}u = (J_{\alpha}^*)^{-1}u
\quad \mbox{if $u \in \WWW$}.
$$

We can generalize the domain of $\pppa$ to $\, _{\beta,p}W(0,T)^*$ with 
arbitrary $\beta > 0$ and refer to Yamamoto \cite{Ya21} for the case
of $p=2$.  Furthermore we can prove several fundamental formulae 
in the fractional calculus such as $\pppa \ppp_t^{\beta} 
= \ppp_t^{\alpha+\beta}$ in $\,_{\alpha+\beta,p}W(0,T)^*$ and 
$W_{\alpha+\beta,p}(0,T)$.  However, we here omit the details and proceed  
to time-fractional ordinary differential equations.
\section{Initial value problems for fractional ordinary differential equations}

We have constructed $\pppa$ in $L^p(0,T)$-based Sobolev spaces in Sections 2 
and 3.  Although we should widely pursue the fractional calculus
of $\pppa$, our main purpose is to study initial value problems for
time-fractional ordinary differential equations and initial boundary value 
problems for time-fractional partial differential equations within 
$L^p$-spaces.  See \cite{Ya21} as for some fractional calculus of $\pppa$ for 
$p=2$.

In this article, we are restricted to initial value problems for
simple fractional ordinary differential equations and illustrate the 
well-posedness of solution within the framework of $\pppa$ in 
$\WWW$ or $\WWWR^*$.

We can make comprehensive studies for general classes, but we consider only 
simple types of single linear fractional ordinary differential equations in 
the case of $0<\alpha<1$.  See Yamamoto \cite{Ya20} as for 
similar treatments
in the case of $p=2$.  One can refer to Diethelm \cite{Die}, 
Jin \cite{J}, which treat various topics concerning fractional derivatives 
and differential equations from different viewpoints.

Throughout this section, we assume that
$$
1 \le p < \infty, \quad 0<\alpha < 1, \quad a\in \R.
$$
\\

{\bf Section 4.1. Single linear fractional differential equation with bounded 
coefficient}

We consider 
$$
\pppa (u(t) - a) = b(t)u + f(t), \quad 0<t<T,   \eqno{(4.1)}
$$
and
$$
u-a \in \WWW.                \eqno{(4.2)}
$$
By the definition of $\pppa$ in $\WWW$, we see that the initial value problem
(4.1) - (4.2) is equivalent to 
$$
u-a = J^{\alpha}(b(t)u) + J^{\alpha}f(t), \quad 0<t<T.      \eqno{(4.3)}
$$
Similarly to \cite{KRY} and \cite{Ya21}, we interpret the initial condition 
by (4.2).  If $\alpha p > 1$, then the Sobolev embedding (e.g.,
\cite{Gri}) and Theorem 2.1 (i) yield that $\WWW \subset 
\, _{0}W^{\alpha-\ep,p}(0,T)\subset C[0,T]$, so that we conclude that 
$u \in C[0,T]$ and $\lim_{t\to 0} u(t) = a$, that is, the initial condition 
(4.2) can be understood in a usual sense.  If $\alpha p \le 1$, then such 
an intepretation is impossible, but we can prove the unique existence 
of solution to (4.1)-(4.2) for all $0<\alpha<1$ and $1\le p < \infty$.
\\
{\bf Theorem 4.1.}
\\
{\it
Let $b\in L^{\infty}(0,T)$ and $f\in L^p(0,T)$.  Then there exists a unique 
solution $u=u(t)$ to (4.1) - (4.2) and we can find a constant $C>$ such that 
$$
\Vert u-a\Vert_{\WWW} \le C(\Vert f\Vert_{\LLLP} + \vert a\vert)
                                       \eqno{(4.4)}
$$
for all $f \in \LLLP$ and $a\in \R$.
}

For $p\ne 2$, we can not completely characterize $\WWW$, 
but Theorem 2 (i) and (ii) provide the properties of the domain $\WWW$ of 
$\pppa$ intermediated between $\, _{0}W^{\alpha+\ep,p}(0,T)$ and 
$\, _{0}W^{\alpha-\ep,p}(0,T)$ with any small gap $\ep>0$, and 
we can apply $\pppa$ for the regularity of solutions
to fractional equations in a flexible way.

If $\alpha p < 1$, then  
$$
a = \frac{a}{\Gamma(\alpha)\Gamma(1-\alpha)}\int^t_0 (t-s)^{\alpha-1}
s^{-\alpha} ds \quad \mbox{for $0<t<T$ and $s^{-\alpha} \in L^p(0,T)$}, 
$$
so that $a \in \WWW$.  Therefore for $\alpha p < 1$, 
the estimate (4.4) is rewritten as 
$$
\Vert u\Vert_{\WWW} \le C(\Vert f\Vert_{\LLLP} + \vert a\vert).
$$
However, if $\alpha p \ge 1$, then we can not obtain the above 
estimate.

In the case $b\in L^{\infty}(0,T)$, Theorem 4.1 can be found in 
the existing works (e.g., Theorem 3.2 (p.124) in \cite{Po}).
We here present Theorem 4.1 in order to illustrate our arguments
which can work for more general cases, 
as are discussed in Sections 4.3 and 4.4 later.
\\
\vspace{0.1cm}
\\
{\bf Proof of Theorem 4.1.}
\\
It suffices to consider (4.3).  Setting $v:= u-a$ and defining an operator
$K: \LLLP \, \RRRR \, \LLLP$ by  
$Kv(t):= J^{\alpha}(b(t)v)$, we see that (4.3) is 
equivalent to 
$$
v(t) = Kv(t) + J^{\alpha}(b(t)a + f(t)), \quad 0<t<T.       \eqno{(4.5)}
$$
The theorem will be proved if we can show that 
the operator $K$ possesses a unique fixed point.

By $b\in L^{\infty}(0,T)$, we deduce that $v\in \LLLP \RRRR bv\in \LLLP$ is 
a bounded operator.  By (2.17), it follows that $v\in \LLLP \, \RRRR\,
J^{\alpha}(bv) \in \LLLP$ is a compact operator.
Therefore, $K: \LLLP \, \RRRR \, \LLLP$ is compact.

Assume that $a=0$ and $f=0$ in (4.5), that is, $v\in \LLLP$ satisfies
$v=Kv$ in $(0,T)$.  Then
$$
v(t) = \frac{1}{\Gamma(\alpha)}\int^t_0 (t-s)^{\alpha-1}v(s) ds,
\quad 0<t<T.
$$
Hence,
$$
\vert v(t)\vert \le C\int^t_0 (t-s)^{\alpha-1}\vert v(s)\vert ds,
\quad 0<t<T.
$$
By the generalized Gronwall inequality (e.g., Lemma 7.1.1 (p.188) in 
Henry \cite{H} or Lemma A.2 in \cite{KRY}), we obtain that $v=0$ in 
$(0,T)$.  Consequently the Fredholm alternative implies that there exists 
a unique fixed point to (4.5) and
$$
\Vert v\Vert_{\LLLP} \le C(\Vert J^{\alpha}(ba+f)\Vert_{\LLLP}
\le C(\vert a\vert + \Vert f\Vert_{\LLLP})
$$
by $b\in L^{\infty}(0,T)$.  Hence, since 
$v=J^{\alpha}(bv+ba+f)$ by (4.5), using the norm of $\Vert \cdot\Vert
_{\WWW}$ and $b\in L^{\infty}(0,T)$, we obtain
$$
\Vert v\Vert_{\WWW} = \Vert bv + ba + f\Vert_{\LLLP}
\le C(\vert a\vert + \Vert f\Vert_{\LLLP}).
$$
Thus the proof of Theorem 4.1 is complete.
$\BLS$

{\bf Section 4.2. Single linear multi-term fractional differential equation 
with bounded coefficient}

Let 
$$
0<\alpha_1 < \cdots < \alpha_N < \alpha < 1.
$$
We consider an initial value problem
$$
\left\{ \begin{array}{rl}
& \pppa (u-a) + \sum_{k=1}^N b_k(t)\ppp_t^{\alpha_k}(u-a)
= b(t)u(t) + f(t), \quad 0<t<T, \\
& u-a \in \WWW.
\end{array}\right.
                                          \eqno{(4.6)}
$$
By Proposition 2.1 (i), we note that $\WWW \subset W_{\alpha_k,p}(0,T)$, and so 
$\ppp_t^{\alpha_k}(u-a)\in \LLLP$ with $k=1,..., N$, are well-defined if
$u-a \in \WWW$.  Now we can prove
\\
\\
{\bf Theorem 4.2.}
\\
{\it
Let $b_1, ..., b_N, b\in L^{\infty}(0,T)$ and $f\in L^p(0,T)$.  Then there exists a unique solution $u=u(t)$ to (4.6) and we can find a constant $C>$ such that the estimate (4.4) holds for all $f \in \LLLP$ and $a\in \R$.
}
\\
{\bf Proof of Theorem 4.2.}
\\
Setting $w:= \pppa (u-a) = (J^{\alpha})^{-1}(u-a)$, by Lemma 2.1 (i) we have
$u-a = J^{\alpha}w$ and
$$
\ppp_t^{\alpha_k}(u-a) = \ppp_t^{\alpha_k}J^{\alpha}w
= \ppp_t^{\alpha_k}(J^{\alpha_k}J^{\alpha-\alpha_k})w 
= J^{\alpha-\alpha_k}w.
$$
Then (4.6) is equivalent to
$$
w(t) = -\sum_{k=1}^N b_k(t)J^{\alpha-\alpha_k}w + b(t)J^{\alpha}w
+ b(t)a + f(t), \quad 0<t<T.
$$
By (2.17) and $b, b_1, ...., b_N \in L^{\infty}(0,T)$, we can verify that 
the operator
$$
Kw(t): = -\sum_{k=1}^N b_k(t)J^{\alpha-\alpha_k}w + b(t)J^{\alpha}w
$$ 
is compact from $\LLLP$ to itself.  Similarly to Theorem 4.1, we can verify
the unique existence of the fixed point to $w=Kw + (ba + f)$.
Thus we can complete the proof of Theorem 4.2.
$\BLS$

{\bf Section 4.3. Single linear fractional differential equation with 
unbounded coefficient}

We return to a simple equation:
$$
\pppa (u-a) = b(t)u, \quad 0<t<T, \qquad 
u-a \in \WWW.                                  \eqno{(4.7)}
$$

Here we study the unique existence and the regularity of the solution to
(4.7) for $b \in L^p(0,T)$ with $1\le p < \infty$ within our framework.

As is discussed in Section 4.1, the regularity $b\in L^{\infty}(0,T)$
makes the total arguments very simple, and for $b\not\in L^{\infty}(0,T)$,
we need more careful discussions.

We remark that the case $\alpha=1$ does not require us any 
special consideration because a formula of solution 
$$
u(t) = a\exp\left(\int^t_0 b(s) ds\right), \quad 0<t<T
$$
implies the unqiue existence of the solution in 
$W^{1,\infty}(0,T)$ for arbitrary $b\in L^1(0,T)$.
However, for $0<\alpha<1$, we do not have such an explicit formula of
solution to (4.7) and we need proper arguments even though the equation 
in (4.7) is extremely simple.
In order to concentrate on the regularity issue of the coefficient 
$b(t)$, we consider (4.7) without the right-hand side of the 
equation.

In terms of $\pppa$ with $\DDD(\pppa) = \WWW$, we prove
\\
{\bf Theorem 4.3.}
\\
{\it
Let $0<\alpha<1$.  We assume that 
$$
b \in \LLLQ  \quad \mbox{with $1<q < \infty$}.       \eqno{(4.8)}
$$
If
$$
1 < p < \infty, \quad \frac{1}{\alpha} < q < \infty, \quad 
p-1 > \frac{1}{q-1},                                 \eqno{(4.9)}
$$
then there exists a unique solution to (4.7) such that 
$$
u-a \in W_{\alpha, \, \frac{pq}{p+q}}(0,T), \quad
u \in \LLLP.                                 \eqno{(4.10)}
$$
Moreover there exists a constant $C>0$ such that 
$$
\Vert u-a\Vert_{W_{\alpha, \, \frac{pq}{p+q}}(0,T)}
+ \Vert u\Vert_{\LLLP} \le C\vert a\vert
$$
for all $a\in \R$.
}

As is seen by the proof below, if $1<q \le \frac{1}{\alpha}$, then our 
proof does not work, and we do not know the unique existence of $u$.
\\
{\bf Proof of Theorem 4.3.}
\\
Setting $v:= u-a$, we see that (4.7) is equivalent to
$$
v = Kv(t) + aJ^{\alpha}b(t), \quad 0<t<T,      \eqno{(4.11)}
$$
where we set $Kv(t):= J^{\alpha}(b(t)v(t))$, $0<t<T$ for 
$v \in \LLLP$.
Theorem 2.1 yields
$$
J^{\alpha}b\in W_{\alpha,q}(0,T) \subset W^{\alpha-\ep,q}(0,T).   
                  \eqno{(4.12)}
$$
Here $\ep>0$ is a sufficiently small constant.

We note that $p-1 > \frac{1}{q-1}$ in (4.9) implies $p>1$.
We set
$$
r = \frac{pq}{p+q}.
$$
Then by the same condition in (4.9), we deduce 
$$
1 < r \le p, q, \quad \frac{1}{r} = \frac{1}{p} + \frac{1}{q}.
                                                   \eqno{(4.13)}
$$
Therefore, for $b \in \LLLQ$, the H\"older inequality implies
$bv \in L^r(0,T)$ for $v \in L^p(0,T)$, Theorem 2.1 (i) implies
$$
J^{\alpha}(bv) \in W^{\alpha-\ep,r}(0,T) \quad 
\mbox{if $v \in L^p(0,T)$}.                     \eqno{(4.14)}
$$
We choose a small constant $\ep > 0$ such that $0 < \ep < \frac{\alpha}{2}$.
Since $q \ge r$ in (4.13) yields $W^{\alpha-\ep,q}(0,T) 
\subset W^{\alpha-\ep,r}(0,T)$, so that by (4.12) we have
$J^{\alpha}b \in W^{\alpha-\ep,r}(0,T)$.

Therefore, by (4.14), we see that if 
$$
W^{\alpha-2\ep,r}(0,T) \subset L^p(0,T),        \eqno{(4.15)}
$$
then 
$$
Kv + aJ^{\alpha}b \in \LLLP \quad \mbox{for each $v \in \LLLP$}.
$$
Moreover, Theorem 2.1 (i) and (4.14) yield that 
$v \in \LLLP\, \RRRR\, J^{\alpha}(bv) \in W^{\alpha-\ep,r}(0,T)$ is
a bounded operator.  Since the embedding $W^{\alpha-\ep,r}(0,T)\, 
\RRRR\, W^{\alpha-2\ep,r}(0,T)$ is compact, it follows that 
under (4.15), the operator $v \in \LLLP \, \RRRR\, Kv \in \LLLP$ is compact.

Thus in view of the Fredholm alternative, similarly to the proof of
Theorem 4.1, it is sufficient to verify that (4.9) implies (4.15).
We consider the following three cases separately.
\\
{\bf Case 1: $\alpha r < 1$:}
\\
Under condition $q > \frac{1}{\alpha}$, we see that $\alpha r <  1$ is 
equivalent to 
$$
p < \frac{q}{\alpha q-1}.                          \eqno{(4.16)}
$$
By $(\alpha-2\ep)r < 1$, the Sobolev embedding implies
$$
W^{\alpha-2\ep,r}(0,T)\,\subset 
L^{\frac{r}{1-(\alpha-2\ep)r}}(0,T).
$$
We can directly verify that $q > \frac{1}{\alpha}$ yields 
$\frac{r}{1-\alpha r} > p$.  Hence, with sufficiently small
$\ep > 0$, we have
$\frac{r}{1-(\alpha-2\ep)r} > p$, that is,
$$
\LLLP\, \supset \, L^{\frac{r}{1-(\alpha-2\ep)r}}(0,T).
$$
Thus (4.15) holds if (4.9) and (4.16) are satisfied.
\\
{\bf Case 2: $\alpha r = 1$.}
\\
We note that $\alpha r = 1$ is equivalent to 
$p = \frac{q}{\alpha q-1}$.
For sufficiently small $\ep > 0$, we have $(\alpha-2\ep)r < 1$ and 
$$
W^{\alpha-2\ep,r}(0,T)\, \supset\, L^{\frac{r}{1-(\alpha-2\ep)r}}(0,T)
$$
by the Sobolev embedding.  Therefore, in the same was as Case 1, we can 
verify (4.15) under (4.9) and
$$
p = \frac{q}{\alpha q-1}.                    \eqno{(4.17)}
$$
\\
{\bf Case 3: $\alpha r > 1$.}
\\
We note that $\alpha r > 1$ is equivalent to 
$p > \frac{q}{\alpha q -1}$.  
Choosing $\ep>0$ sufficiently small, we obtain $(\alpha-2\ep)r > 1$.
Therefore, the Sobolev embedding yields that $W^{\alpha-2\ep,r}(0,T)
\subset L^{\infty}(0,T)$.  Hence, (4.15) holds if (4.9) and 
$$
p > \frac{q}{\alpha q -1}              \eqno{(4.18)}
$$
are satisfied.
\\

Thus, taking the union of the sets of $(p,q)$ satisfying (4.16) - (4.18), 
we verify that (4.15) holds if (4.9) is satisfied.
Thus we complete the proof of Theorem 4.3.
$\BLS$
 
{\bf Section 4.4. Single linear fractional differential equation with 
singular non-homogeneous term}

Let $0<\alpha<1$, $\beta \ge 0$ and $1 <  p \le \infty$.  Henceforth we define 
$1 \le q < \infty$ by 
$$
p= 
\left\{ \begin{array}{rl}
& \frac{q}{q-1} \quad \mbox{if $1<q<\infty$}, \\
&\infty \quad \mbox{if $q=1$}.
\end{array}\right.
$$
We consider 
$$
\left\{ \begin{array}{rl}
& \pppa (u-a) = b(t)u(t) + f(t), \quad 0<t<T, \\
& u-a \in L^p(0,T),
\end{array}\right.
                                         \eqno{(4.19)}
$$
where $b\in L^{\infty}(0,T)$ and 
$f \in \, _{\alpha,q}W(0,T)^*$.

In (4.19), we understand as 
$$
\pppa = (J_{\alpha}^*)^{-1}:\, L^q(0,T)^* \, \RRRR\,
 \, _{\alpha,q}W(0,T)^*,
$$
that is,
$$
\left\{ \begin{array}{rl}
& \pppa :\, \LLLP \, \RRRR\, \, _{\alpha,q}W(0,T)^* \quad
\mbox{if $1<p<\infty$}, \\
& \pppa :\, L^{\infty}(0,T) \, \RRRR\, \, _{\alpha,1}W(0,T)^* \quad
\mbox{if $p=\infty$}.
\end{array}\right.
                                            \eqno{(4.20)}
$$
Moreover, by Lemma 3.4 (iii), we see that 
$$
J_{\alpha}^*:\, _{\alpha+\beta,q}W(0,T)^* \, \RRRR \, 
_{\beta,q}W(0,T)^*
$$
is an isomorphism.

In terms of (4.20), by setting $v:= u-a$, the initial value problem
(4.19) is equivalent to
$$
\left\{ \begin{array}{rl}
& v = J_{\alpha}^*(bv) + J_{\alpha}^*(ba+f), \\
& v \in \LLLP.
\end{array}\right.
                                     \eqno{(4.21)}
$$
For simplicity, we assume that $1<p<\infty$.
Now we can prove
\\
{\bf Theorem 4.4.}
\\
{\it 
Let $b\in L^{\infty}(0,T)$ be arbitrary.  Then there exists a unique 
solution $u$ to (4.19) and we can find a constant $C>0$ such that  
$$
\Vert u-a\Vert_{\LLLP} \le C(\vert a\vert 
+ \Vert f\Vert_{_{\alpha,q}W(0,T)^*})
$$
for all $a \in \R$ and $f \in \,_{\alpha,q}W(0,T)^*$.
}
\\
{\bf Example.}
\\
Given constant $t_0\in [0,T]$, we consider a Dirac delta function 
$f(t):= \delta_{t_0}$, that is, \\
$\, _{C[0,T]^*}< \delta_{t_0},\, \psi>_{C[0,T]} 
:= \psi(t_0)$ for all $\psi \in C[0,T]$.
Let $\alpha q > 1$ with $q = \frac{p}{p-1}$.  Then the Sobolev
embedding yields that $\, _{\alpha,q}W(0,T) \subset C[0,T]$, and so
it turns that 
$\delta_{t_0} \in \, W_{\alpha,q}(0,T)^*$.
Therefore Theorem 4.4 asserts that there exists a unique solution 
$u\in \LLLP$ to $\pppa (u-a) = b(t)u(t) + \delta_{t_0}(t)$ 
in $\, _{\alpha,q}W(0,T)^*$.
\\
\vspace{0.1cm}
\\
{\bf Proof of Theorem 4.4.}
\\
By Lemma 3.4 (iii), noting that 
$\, _{0,q}W(0,T)^* = L^q(0,T)^* = L^p(0,T)$, we obtain
$$
\Vert J_{\alpha}^*f\Vert_{\LLLP}
\le C\Vert f\Vert_{_{\alpha,q}W(0,T)^*}.
$$
By $b \in L^{\infty}(0,T)$, we apply (3.8) to have
$J_{\alpha}^*b = J^{\alpha}b \in L^{\infty}(0,T)$ by the Young 
inequality.  Therefore,
$$
J_{\alpha}^*(ba+f) \in \LLLP, \quad
\Vert J_{\alpha}^*(ba+f)\Vert_{\LLLP}
\le C(\vert a\vert + \Vert f\Vert_{_{\alpha,q}W(0,T)^*}).
                                                         \eqno{(4.22)}
$$
Since $b\in L^{\infty}(0,T)$ implies $bv \in \LLLP$.  Again by
(3.8), we deduce that $J_{\alpha}^*(bv) = J^{\alpha}(bv)$ for all
$v\in \LLLP$.  Thus (4.21) is rewritten as 
$$
v = J^{\alpha}(bv) + G, \quad v\in \LLLP,
$$
where $G:= J^{\alpha}(ba+f) \in \LLLP$. 
Consequently our argument can be executed within $\LLLP$, and we can 
repeat the proof of Theorem 4.1.
Thus the proof of Theorem 4.4 is complete.
$\BLS$
\section{Concluding remarks}

{\bf 1.} In Section 2 of this article, we have established a time-fractional 
derivative $\pppa$ in $L^p(0,T)$-based Sobolev-Slobodecki spaces for 
$1\le p < \infty$.  Theorem 2.2 provides characterization of
the domain $\DDD(\pppa)$ in terms of the Sobolev-Slobodecki spaces,
which still admits a gap with any small $\ep>0$ in Sobolev orders.

We remark that the domain $\DDD(\pppa)$ in a special case $p=2$, is
completely characterized by the Sobolev-Slobodecki spaces 
(see \cite{GLY}, \cite{KRY}, \cite{Ya21}).

{\bf 2.} By the duality, in Section 3, for $1\le p < \infty$,
we extend the domain $\DDD(\pppa)$ to $\LLLP$ with the range in the dual space 
$\WWWR^*$.

{\bf 3.}  In Section 4, we established the unique existence of solutions to 
initial value problems for fractional ordinary differential equations on the 
basis of $\pppa$.  We can develop more comprehensive treatments for wider
classes of equations, but we postpone them and are restricted to single
fractional equations.  Even for such simple equations, there are no works
for example in the case where a coefficient $b(t)$ is not bounded or 
non-homogeneous term $f(t)$ does not belong to $\LLLP$, 
and we established the unique
existences of solution for $b \in L^p(0,T)$ with $p\ne \infty$ or
$f \in \, _{\alpha,q}W(0,T)^*$ with $\alpha > 0$.

Furthermore we here considered $\pppa$ only in the case of $0<\alpha<1$ 
although in Sections 2 and 3 the fractional derivative $\pppa$ is defined for 
all $\alpha > 0, \not\in \N$, and more general treatments will be provided in 
a future work.

We recall that Theorems 2.1 and 2.2 give characterization of the domain 
of $\pppa$ which is not the best possible because of regularity loss $\ep>0$, 
even though $\ep>0$ can be choosen  arbitrarily small.  
However, for studying fractional
differential equations, such characterization of the domain by 
the Sobolev spaces is used mainly for applying the Sobolev embedding, and, 
as is illustrated in Section 4, we emphasize that the $\ep$-gap is 
not a serious disadvantage.

{\bf 4.}  We do not discuss time-fractional partial differential equations 
at all.  The related topics should be studied also in future works
for $p \ge 1$.  We remark that the case $p=2$ is relatively satisfactorily 
argued already in e.g., \cite{GLY}, \cite{KRY}.
\section{Appendix: Proof of Lemma 2.4.}

For $m=1$ and $p=1$, Lemma 2.4 is proved as 
Theorem 4.2.3 (p.77) in \cite{GV}, and our
proof is based on the expression (6.1) (which is used also on p.77
of \cite{GV}).
\\
\vspace{0.1cm}
\\
{\bf First Step: $0<\alpha<1$.}
\\
We note that $\CC \subset W^{\alpha+\ep,p}(0,T)$ if 
$0 < \alpha + \ep < 1$.
Let $u \in \CC$.  Then
\begin{align*}
& J^{1-\alpha}u(t) = \frac{1}{\Gamma(1-\alpha)}
\int^t_0 (t-s)^{-\alpha} u(s) ds\\
=& \frac{1}{\Gamma(1-\alpha)}\int^t_0 (t-s)^{-\alpha}
(u(s)-u(t)) ds + \frac{u(t)}{\Gamma(1-\alpha)}\int^t_0 (t-s)^{-\alpha}ds,
\end{align*}
that is,
$$
J^{1-\alpha}u(t) 
= \frac{1}{\Gamma(1-\alpha)}\int^t_0 (t-s)^{-\alpha}
(u(s)-u(t)) ds + \frac{u(t)t^{1-\alpha}}{(1-\alpha)\Gamma(1-\alpha)},
\quad 0<t<T.                                  
$$
Since $u \in \CC$, we see
$$
\vert (t-s)^{-\alpha}(u(s) - u(t))\vert \le C\vert t-s\vert^{1-\alpha}
$$
and 
\begin{align*}
& \left\vert \frac{\ppp}{\ppp t}((t-s)^{-\alpha}(u(s) - u(t)))
\right\vert\\ 
= & \left\vert -\alpha (t-s)^{-\alpha-1}(u(s) - u(t))
- (t-s)^{-\alpha}\frac{du}{dt}(t)\right\vert
\le C\vert t-s\vert^{-\alpha}
\end{align*}
so that $(t-s)^{-\alpha}(u(s) - u(t))$ is in $W^{1,1}(0,t)$ as a function
in $s \in (0,t)$ for arbitrarily fixed $t \in (0,T)$.
Therefore,
\begin{align*}
& \frac{\ppp}{\ppp t}\int^t_0 (t-s)^{-\alpha}(u(s) - u(t)) ds\\
=& -\alpha \int^t_0 (t-s)^{-\alpha-1}(u(s) - u(t)) ds
+ \int^t_0 (t-s)^{-\alpha}\left( -\frac{du}{dt}(t)\right) ds \\
=& -\alpha \int^t_0 (t-s)^{-\alpha-1}(u(s) - u(t)) ds
- \frac{t^{1-\alpha}}{1-\alpha}\frac{du}{dt}(t).
\end{align*}
Hence,
$$
\frac{d}{dt}J^{1-\alpha}u(t) 
= -\alpha \int^t_0 (t-s)^{-\alpha-1}(u(s) - u(t)) ds
+ \frac{t^{-\alpha}}{\Gamma(1-\alpha)}u(t)
$$
$$
=: S_1(t) + S_2(t), \quad 0<t<T  \quad
\mbox{for $u \in \CC$.}                   \eqno{(6.1)}
$$
\\
{\bf Second Step: Estimation of $S_1(t)$ for $0<\alpha<1$ and
$u\in \CC$.}
\\
{\bf Case 1: $p=1$.}
\\
We have
$$
\Vert S_1\Vert_{L^1(0,T)} \le 
\int^T_0 \left\vert \int^t_0 (t-s)^{-\alpha-1}
(u(s) - u(t)) ds \right\vert dt
$$
$$
\le \int^T_0\int^T_0 \frac{\vert u(t) - u(s)\vert}
{\vert t-s\vert^{\alpha+1}} dsdt \le \Vert u\Vert_{W^{\alpha,1}(0,T)}.
                                       \eqno{(6.2)}
$$
\\
{\bf Case 2: $1<p<\infty$.}
\\
We have
$$
\left\vert \frac{\vert u(t) - u(s)\vert}{\vert t-s\vert^{\alpha+1}}
\right\vert
= \frac{\vert u(t) - u(s)\vert}{\vert t-s\vert^{\frac{1}{p}+\alpha+\ep}}
\frac{1}{\vert t-s\vert^{1-\frac{1}{p}-\ep}}.
$$
Set $\frac{1}{q}:= 1 - \frac{1}{p}$, that is, 
$q = \frac{p}{p-1} > 1$.  Then 
$$
\int^T_0 \vert t-s\vert^{q\left(\frac{1}{p}+\ep-1\right)} ds
= \int^T_0 \vert t-s\vert^{-1+\frac{p\ep}{p-1}} ds
\le \left( T^{\frac{p\ep}{p-1}} + (T-t)^{\frac{p\ep}{p-1}}\right)
\frac{p-1}{p\ep} < \infty.
$$
Consequently, the H\"older inequality yields 
\begin{align*}
& \int^T_0 \frac{\vert u(t) - u(s)\vert}{\vert t-s\vert^{\alpha+1}} ds
= \int^T_0 \frac{\vert u(t) - u(s)\vert}
{\vert t-s\vert^{\frac{1}{p}+\alpha+\ep}} 
\frac{1}{\vert t-s\vert^{1-\frac{1}{p}-\ep}} ds\\
\le& \left( \int^T_0 \frac{\vert u(t) - u(s)\vert^p}
{\vert t-s\vert^{1+p(\alpha+\ep)}} ds\right)^{\frac{1}{p}}
\left( \int^T_0 \vert t-s\vert^{q\left(\frac{1}{p}+\ep-1\right)} ds
\right)^{\frac{1}{q}}
\le C\left( \int^T_0 \frac{\vert u(t) - u(s)\vert^p}
{\vert t-s\vert^{1+p(\alpha+\ep)}} ds\right)^{\frac{1}{p}}.
\end{align*}
Therefore,
$$
\int^T_0 \vert S_1(t)\vert^p dt 
\le C\int^T_0 \left\vert \int^t_0 \frac{\vert u(t) - u(s)\vert}
{\vert t-s\vert^{1+\alpha}} ds \right\vert^p dt
$$
$$
\le C\int^T_0 \int^T_0 \frac{\vert u(t) - u(s)\vert^p}
{\vert t-s\vert^{1+p(\alpha+\ep)}} dsdt 
\le C\Vert u\Vert_{W^{\alpha+\ep,p}(0,T)}.                      \eqno{(6.3)}
$$
By (6.2) and (6.3), for $p \ge 1$, we obtain
$$
\Vert S_1\Vert_{\LLLP} \le C\Vert u\Vert_{W^{\alpha+\ep,p}(0,T)}.
                                                                \eqno{(6.4)}
$$
\\
{\bf Third Step: Estimation of $S_2(t)$ for $0<\alpha<1$ and
$u \in \CC$.}
\\
We have
$$
\Vert S_2\Vert^p_{\LLLP} \le C\int^T_0 t^{-p\alpha}\vert u(t)\vert^p dt.
$$
\\
{\bf Case 1: $p(\alpha+\ep) < 1$}.
\\
We can have the following Sobolev embedding: 
$$
W^{\alpha+\ep, p}(0,T) \subset  
L^{\frac{p(1-\delta_0)}{1-p(\alpha+\ep)}}(0,T),       \eqno{(6.5)}
$$
where we fix a sufficiently small constant $\delta_0>0$.
The inclusion (6.5) is proved in Adams \cite{Ad} or Grisvard \cite{Gri}
for $p>1$ and as Lemma 4.2.1 (pp.72-73) in Gorenflo and Vessella \cite{GV}
for $p=1$.

We set $q:= \frac{1-\delta_0}{1-p(\alpha+\ep)}$ and 
$r:= \frac{1-\delta_0}{p(\alpha+\ep)-\delta_0}$. 
By $p(\alpha+\ep) < 1$, choosing sufficiently 
small $\delta_0 > 0$, see that $q, r > 1$ 
and $\frac{1}{q} + \frac{1}{r} = 1$. 
Hence the H\"older inequality yields
\begin{align*}
& \int^T_0 t^{-\alpha p} \vert u(t)\vert^p dt 
\le \left( \int^T_0 \vert u(t)\vert^{pq} dt \right)^{\frac{1}{q}}
 \left( \int^T_0 t^{-\alpha pr} dt \right)^{\frac{1}{r}}\\
\le& \left( \int^T_0 \vert u(t)\vert^{\frac{p(1-\delta_0)}{1-p(\alpha+\ep)}} dt
\right)^{\frac{1}{q}}
\left( \int^T_0 t^{-\alpha p\frac{1-\delta_0}{p(\alpha+\ep)-\delta_0}} dt 
\right)^{\frac{1}{r}}.
\end{align*}
In view of $\alpha p < p(\alpha+\ep) < 1$, choosing $\delta_0 > 0$ smaller,
we can assume that 
$\delta_0 < \frac{p\ep}{1-\alpha p}$.  Then direct calculations imply
$$
\alpha p \frac{1-\delta_0}{p(\alpha+\ep)-\delta_0} < 1,
$$
and so 
$$
\left( \int^T_0 t^{ -\alpha p \frac{1-\delta_0}{p(\alpha+\ep)-\delta_0} }dt
\right)^{\frac{1}{r}} < \infty.
$$
With (6.5), we obtain
$$
\Vert S_2\Vert^p_{\LLLP} 
\le C\left( \Vert u\Vert_{W^{\alpha+\ep,p}(0,T)}
^{\frac{p(1-\delta_0)}{1-p(\alpha+\ep)}} \right)^{\frac{1}{q}}
= C\Vert u\Vert_{W^{\alpha+\ep,p}(0,T)}^p.                \eqno{(6.6)}
$$
\\
{\bf Case 2: $p(\alpha+\ep)=1$}.
\\
By $\alpha+\ep < 1$, we have $p>1$.  Then the Sobolev embedding 
(e.g., \cite{Ad}, \cite{Gri}) implies
$$
W^{\alpha+\ep,p}(0,T) \subset L^{pq}(0,T)             \eqno{(6.7)}
$$
for any $q\ge 1$.  We set $r := \frac{q}{q-1}$ with $q>1$.
Choosing $q>1$ sufficiently large, since $\lim_{q\to\infty} r = 1$ and 
$\alpha p < \alpha(p+\ep) = 1$, we can obtain $\alpha pr < 1$.
Therefore, (6.7) and the H\"older inequality imply 
\begin{align*}
& \int^T_0 t^{-\alpha p} \vert u(t)\vert^p dt 
\le \left( \int^T_0 \vert u(t)\vert^{pq} dt \right)^{\frac{1}{q}}
 \left( \int^T_0 t^{-\alpha pr} dt \right)^{\frac{1}{r}}\\
\le& C\Vert u\Vert^p_{L^{pq}(0,T)} 
\le C\Vert u\Vert^p_{W^{\alpha+\ep,p}(0,T)},
\end{align*}
that is, we reach (6.6) in Case 2.
\\
{\bf Case 3: $p(\alpha+\ep) > 1$}.
\\
Then we note that $p>1$.  The Sobolev embedding yields
$$
W^{\alpha+\ep,p}(0,T) \subset C^{\theta}[0,T],      \eqno{(6.8)}
$$
where
$$
0 < \theta < \alpha+\ep - \frac{1}{p}
$$
(e.g., \cite{Gri}).
Here we set 
\begin{align*}
&C^{\theta}[0,T] := \{ v\in C[0,T];\, 
\mbox{there exists a constant $C=C_v > 0$ such that}\\
& \vert v(t) - v(s) \vert \le C_v\vert t-s\vert^{\theta}\}
\end{align*}
and
$$
\Vert v\Vert_{C^{\theta}[0,T]}:= \Vert v\Vert_{C[0,T]}
+ \sup_{t,s\in [0,T], t\ne s} \frac{\vert u(t)-u(s)\vert}
{\vert t-s\vert^{\theta}}.
$$
By $u \in \CC$, we have $u(0) = 0$ and
$$
\vert u(t)\vert = 
\vert u(t) - u(0)\vert \le \Vert u\Vert_{C^{\theta}[0,T]}t^{\theta}.
$$
Therefore, 
$$
\Vert S_2\Vert^p_{\LLLP} \le C\int^T_0 \frac{\vert u(t)\vert^p}
{t^{\alpha p} } dt 
\le C\int^T_0 t^{p(\theta-\alpha)}\Vert u\Vert_{C^{\theta}[0,T]}^p dt.
$$
Since $0 < \theta < \alpha+\ep-\frac{1}{p}$, we choose $0 < \delta_1 < \ep$ and
set $\theta := \alpha + \ep - \frac{1}{p} - \delta_1$.
Then $p(\theta-\alpha) = -1+p(\ep-\delta_1) > -1$, and 
$\int^T_0 t^{p(\theta-\alpha)} dt < \infty$.
Consequently, by (6.8) we reach 
$$
\Vert u\Vert^p_{\LLLP} \le C\Vert u\Vert_{W^{\alpha+\ep,p}(0,T)}^p,
$$
and (6.6) is proved for Case 3.
\\

Thus, in view of (6.1), (6.4) and (6.6), we finished the proof of 
$$
\Vert J^{1-\alpha}u\Vert_{W^{1,p}(0,T)} 
\le C\Vert u\Vert_{W^{\alpha+\ep,p}(0,T)} \quad \mbox{for 
$u \in \CC$}                          \eqno{(6.9)}
$$
with $0<\alpha<1$, 
where $p\ge 1$ and the constant $C>0$ depends only on $T$, $p$, $\alpha$,
and $\ep$.
\\
{\bf Fourth Step: Completion of the proof for $0<\alpha<1$.}
\\
Let $u\in \, _{0}W^{\alpha+\ep,p}(0,T)$ be arbitrary.  By the definition, 
we can choose $u_n \in \CC$, $n\in \N$ such that $u_n \RRRR u$ in 
$W^{\alpha+\ep,p}(0,T)$ as $n \to \infty$.
In view of (6.9), we deduce that 
$$
\Vert J^{1-\alpha}u_n - J^{1-\alpha}u_m\Vert_{W^{1,p}(0,T)}
\le C\Vert u_n - u_m\Vert_{W^{\alpha+\ep,p}(0,T)} \,
\RRRR \, 0
$$
as $n,m \to \infty$.  Therefore, again by (6.9) and 
$J^{1-\alpha}u_n \in \CC$ for $n\in \N$, there exists 
$\www{w} \in \, _{0}W^{\alpha+\ep,p}(0,T)$ such that 
$$
J^{1-\alpha}u_n \, \RRRR \, \www{w} \quad \mbox{in
$W^{1,p}(0,T)$ as $n \to \infty$},                  \eqno{(6.10)}
$$
and 
$$
\Vert \www{w}\Vert_{W^{1,p}(0,T)} 
= \lim_{n\to \infty} \Vert J^{1-\alpha}u_n\Vert_{W^{1,p}(0,T)}
\le C\lim_{n\to \infty} \Vert u_n\Vert_{W^{\alpha+\ep,p}(0,T)}.
                                               \eqno{(6.11)}
$$
Here we applied the part of the lemma already proved for $0<\alpha+\ep<1$ and 
$u_n \in \CC$.

Since $u_n \RRRR u$ in $W^{\alpha+\ep,p}(0,T)$, implies
that $u_n\, \RRRR\, u$ in $\LLLP$, by the Young inequality, we obtain that 
$J^{1-\alpha}u_n\, \RRRR\, J^{1-\alpha}u$ in $\LLLP$ as
$n \to \infty$.

Combining with (6.10), we see that $J^{1-\alpha}u = \www{w}$.
Therefore, in terms of (6.11), we reach the conclusion of 
Lemma 2.4 for $0<\alpha<1$ and $u\in \, _{0}W^{\alpha+\ep,p}(0,T)$.  
$\BLS$
\\
{\bf Fifth Step: Completion of the proof of Lemma 2.4.}
\\
Let $\alpha = m-1 + \sigma$ with $m=2,3,4,...$ and $0<\sigma < 1$.
By the density argument, similarly to Fourth Step, it suffices to prove
for $u \in \, _{0}C^m[0,T]$. 

Let $u \in \, _{0}C^m[0,T]$ be arbitrarily given.  Then
$$
\frac{d^{m-1}u}{dt^{m-1}} \in \, \CC
$$
by the definition of $\,_{0}W^{\alpha,p}(0,T)$.

Since the conclusion of Lemma 2.4 with $0<\alpha<1$ and 
$u \in \, \CC$ was already proved, we replace 
$m:=1$ and 
$u:= \frac{d^{m-1}u}{dt^{m-1}} \in \, \CC$ to apply, so that 
we obtain
$$
\left\Vert J^{1-\sigma} \frac{d^{m-1}u}{dt^{m-1}}\right\Vert
_{W^{1,p}(0,T)}
\le C\left\Vert \frac{d^{m-1}u}{dt^{m-1}}\right\Vert
_{W^{\sigma+\ep,p}(0,T)}.                      \eqno{(6.12)}
$$
Similarly to (2.16), we can see
$$
\frac{d^{m-1}}{dt^{m-1}}J^{1-\sigma}u 
= J^{1-\sigma}\frac{d^{m-1}u}{dt^{m-1}}\quad \mbox{for $u \in 
\, _{0}C^m[0,T]$}.
$$
Consequently, (6.12) yields
\begin{align*}
& \left\Vert \frac{d^m}{dt^m}J^{1-\sigma}u\right\Vert_{\LLLP}
\le \left\Vert \frac{d^{m-1}}{dt^{m-1}}J^{1-\sigma}u\right\Vert
_{W^{1,p}(0,T)}\\
= & \left\Vert J^{1-\sigma}\left( \frac{d^{m-1}u}{dt^{m-1}}\right)
\right\Vert_{W^{1,p}(0,T)}
\le C\left\Vert \frac{d^{m-1}u}{dt^{m-1}}\right\Vert
_{W^{\sigma+\ep,p}(0,T)}.
\end{align*}
Hence, the definition of the norm $\Vert \cdot\Vert_{W^{\alpha+\ep,p}(0,T)}$
yields
$$
\left\Vert \frac{d^m}{dt^m}J^{1-\sigma}u\right\Vert_{\LLLP}
\le C\Vert u\Vert_{W^{m-1+\sigma+\ep,p}(0,T)}.
                                                        \eqno{(6.13)}
$$
Since $\frac{d^k}{dt^k}(J^{1-\sigma}u)(0) = 0$
by $\frac{d^ku}{dt^k}(0) = 0$ for $0\le k \le m-1$, we have
$$
\frac{d^{m-k}}{dt^{m-k}}(J^{1-\sigma}u)(t) 
= \int^t_0 \frac{(t-s)^{k-1}}{(k-1)!} 
\frac{d^m}{dt^m}(J^{1-\sigma}u)(s) ds, \quad k=1, ..., m
$$
and so 
$$
\Vert J^{1-\sigma}u\Vert_{W^{m-1,p}(0,T)}
\le C\left\Vert \frac{d^m}{dt^m}J^{1-\sigma}u\right\Vert_{\LLLP}.
$$
In view of (6.13) and $m-\alpha = 1-\sigma$, we reach 
$$
\Vert J^{m-\alpha}u\Vert_{W^{m,p}(0,T)}
= \Vert J^{1-\sigma}u\Vert_{W^{m,p}(0,T)}
\le C\Vert u\Vert_{W^{\alpha+\ep,p}(0,T)}
$$
for $u \in \, _{0}C^m[0,T]$.
Thus the proof of Lemma 2.4 is complete.
$\BLS$
\\

{\bf Acknowledgements}.
\\ 
The author is supported by Grant-in-Aid
for Scientific Research (A) 20H00117, JSPS and by the National Natural
Science Foundation of China (Nos.\! 11771270, 91730303).
This paper has been supported by the RUDN University 
Strategic Academic Leadership Program.


\begin{thebibliography}{99}

\bibitem{Ad}
R.A. Adams, {\it Sobolev Spaces}, Academic Press, New York, 1975.

\bibitem{Bre}
H. Brezis, {\it Functional Analysis, Sobolev Spaces and Partial
Differential Equations}, Springer, Berlin, 2011.

\bibitem{CC}
A. Carbotti and G.E. Comi, A note on Riemann-Liouville fractional 
Sobolev spaces, Comm. Pure Applied. Anal. {\bf 20} (2021) 17-54.

\bibitem{Die}
K. Diethelm, {\it The Analysis of Fractional Differential Equations},
Springer, Berlin, 2010.

\bibitem{GKMR}
R. Gorenflo, A.A. Kilbas, F. Mainardi and S. V. Rogosin, {\it Mittag-Leffler 
Functions, Related Topics and Applications}, Springer-Verlag, Berlin, 2014.

\bibitem{GLY}
R. Gorenflo, Yu. Luchko and M. Yamamoto, Time-fractional diffusion equation 
in the fractional Sobolev spaces, Fract. Calc. Appl. Anal. \textbf{18}
(2015) 799--820. 

\bibitem{GV}
R. Gorenflo and S. Vessella, 
{\it Abel Integral Equations}, Lec. Notes in Math. 1461,
Springer, Berlin, 1991.

\bibitem{Gri}
P. Grisvard, {\it Elliptic Problems in Nonsmooth Domains},
Pitman, Boston, 1985. 

\bibitem{H}
D. Henry, {\it Geometric Theory of Semilinear Parabolic Equations},
Lec. Notes in Math. 840, Springer, Berlin, 1981.

\bibitem{J}
B. Jin, {\it Fractional Differential Equations, 
An Approach via Fractional Derivatives}, Springer, Berlin, 2021.

\bibitem{Ki}
Y. Kian, Equivalence of definitions of solutions for some class of fractional 
diffusion equations, 2021, preprint arXiv:2111.06168

\bibitem{KiYa}
Y. Kian and M. Yamamoto, 
Well-posedness for weak and strong solutions of non-homogeneous 
initial boundary value problems for fractional diffusion equations,
Fract. Calc. Appl. Anal. {\bf 24} (2021) 168-201.

\bibitem{KST}
A.A. Kilbas, H.M. Srivastava and J.J. Trujillo, 
{\it Theory and Applications
of Fractional Differential Equations}. Elsevier, Amsterdam, 2006.

\bibitem{Koe}
H. K\"onig, Grenzordnungen von Operatorenidealen (I), (II), Math. Ann.
{\bf 212} (1974) 51-64, 65-77.

\bibitem{KRY}
A. Kubica, K. Ryszewska and M. Yamamoto,   
{\it Time-fractional Differential Equations
A Theoretical Introduction}, Springer, Tokyo, 2020.

\bibitem{LL1}
L. Li and J-G. Li, 
A generalized definition of Caputo detivatives and its application to
fractional ODES, SIAM J. Math. Anal. {\bf 50} (2018) 2867-2900.

\bibitem{LL2}
L. Li and J-G. Li, 
Some compactness criteria for weak solutions of time fractional 
PDEs, SIAM J. Math. Anal. {\bf 50} (2018) 3963-3995.

\bibitem{LY}
Yu. Luchko and M. Yamamoto,
Comparison principle and monotone method for time-fractional diffusion
equations with Robin boundary condition, 2021, preprint.

\bibitem{MK}
R. Metzler and J. Klafter, The random walk's guide to anomalous
diffusion: a fractional dynamics approach, Phys. Rep. {\bf 339}
(2000) 1-77.

\bibitem{Po}
I. Podlubny, {\it Fractional Differential Equations}, Academic Press,  
San Diego, 1999.

\bibitem{SY}
K. Sakamoto and M. Yamamoto, Initial value/boundary value problems for 
fractional diffusion-wave equations and applications to some inverse problems,
J. Math. Anal. Appl. {\bf 382} (2011) 426-447.

\bibitem{Ya20}
M. Yamamoto, 
On time fractional derivatives in fractional Sobolev 
spaces and applications to fractional ordinary differential equations,   
Nonlocal and fractional operators, pp. 287-308, SEMA SIMAI Springer Ser., 26, 
Springer, Cham, 2021.

\bibitem{Ya21}
M. Yamamoto,
Fractional calculus and time-fractional differential equations: 
revisit and construction of a theory, 2022, preprint.

\bibitem{Za}
R. Zacher, Weak solutions of abstract evolutionary integro-differential
equations in Hilbert spaces, Funkcialaj Ekvacioj {\bf 52} (2009) 1-18.
\end{thebibliography}
\end{document}